%
\catcode`@=11
%
%
\def\bibn@me{R\'ef\'erences}
\def\bibliographym@rk{\centerline{{\sc\bibn@me}}
    \sectionmark\section{\ignorespaces}{\unskip\bibn@me}
    \bigbreak\bgroup
    \ifx\ninepoint\undefined\relax\else\ninepoint\fi}
%
%
%
\let\refsp@ce=\
\let\bibleftm@rk=[
\let\bibrightm@rk=]
%
%
%
\def\numero{n\raise.82ex\hbox{$\fam0\scriptscriptstyle o$}~\ignorespaces}
%
%
\newcount\equationc@unt
\newcount\bibc@unt
\newif\ifref@changes\ref@changesfalse
\newif\ifpageref@changes\ref@changesfalse
\newif\ifbib@changes\bib@changesfalse
\newif\ifref@undefined\ref@undefinedfalse
\newif\ifpageref@undefined\ref@undefinedfalse
\newif\ifbib@undefined\bib@undefinedfalse
\newwrite\@auxout
%
%
\def\eqnum{\global\advance\equationc@unt by 1%
\edef\lastref{\number\equationc@unt}%
\eqno{(\lastref)}}
%
%
%
%
%
%
\def\re@dreferences#1#2{{%
    \re@dreferenceslist{#1}#2,\undefined\@@}}
\def\re@dreferenceslist#1#2,#3\@@{\def\next{#2}%
    \expandafter\ifx\csname#1@@\meaning\next\endcsname\relax
    ??\immediate\write16
    {Warning, #1-reference "\next" on page \the\pageno\space
    is undefined.}%
    \global\csname#1@undefinedtrue\endcsname
    \else\csname#1@@\meaning\next\endcsname\fi
    \ifx#3\undefined\relax
    \else,\refsp@ce\re@dreferenceslist{#1}#3\@@\fi}
%
%
%
\def\newlabel#1#2{{\def\next{#1}\newl@bel#2}}
\def\newl@bel#1#2{%
    \expandafter\xdef\csname ref@@\meaning\next\endcsname{#1}%
    \expandafter\xdef\csname pageref@@\meaning\next\endcsname{#2}}
\def\label#1{{%
    \toks0={#1}\message{ref(\lastref) \the\toks0,}%
    \ignorespaces\immediate\write\@auxout%
    {\noexpand\newlabel{\the\toks0}{{\lastref}{\the\pageno}}}%
    \def\next{#1}%
    \expandafter\ifx\csname ref@@\meaning\next\endcsname\lastref%
    \else\global\ref@changestrue\fi%
    \newlabel{#1}{{\lastref}{\the\pageno}}}}
\def\ref#1{\re@dreferences{ref}{#1}}
\def\pageref#1{\re@dreferences{pageref}{#1}}
%
%
\def\bibcite#1#2{{\def\next{#1}%
    \expandafter\xdef\csname bib@@\meaning\next\endcsname{#2}}}
\def\cite#1{\bibleftm@rk\re@dreferences{bib}{#1}\bibrightm@rk}
%
%
\def\beginthebibliography#1{\bibliographym@rk
    \setbox0\hbox{\bibleftm@rk#1\bibrightm@rk\enspace}
    \parindent=\wd0
    \global\bibc@unt=0
    \def\bibitem##1{\global\advance\bibc@unt by 1
        \edef\lastref{\number\bibc@unt}
        {\toks0={##1}
        \message{bib[\lastref] \the\toks0,}%
        \immediate\write\@auxout
        {\noexpand\bibcite{\the\toks0}{\lastref}}}
        \def\next{##1}%
        \expandafter\ifx
        \csname bib@@\meaning\next\endcsname\lastref
        \else\global\bib@changestrue\fi%
        \bibcite{##1}{\lastref}
        \medbreak
        \item{\hfill\bibleftm@rk\lastref\bibrightm@rk}%
        }
    }
\def\endthebibliography{\egroup\par}
%
%
%
\def\@closeaux{\closeout\@auxout
    \ifref@changes\immediate\write16
    {Warning, changes in references.}\fi
    \ifpageref@changes\immediate\write16
    {Warning, changes in page references.}\fi
    \ifbib@changes\immediate\write16
    {Warning, changes in bibliography.}\fi
    \ifref@undefined\immediate\write16
    {Warning, references undefined.}\fi
    \ifpageref@undefined\immediate\write16
    {Warning, page references undefined.}\fi
    \ifbib@undefined\immediate\write16
    {Warning, citations undefined.}\fi}
%
%
\immediate\openin\@auxout=\jobname.aux
\ifeof\@auxout \immediate\write16
  {Creating file \jobname.aux}
\immediate\closein\@auxout
\immediate\openout\@auxout=\jobname.aux
\immediate\write\@auxout {\relax}%
\immediate\closeout\@auxout
\else\immediate\closein\@auxout\fi
%
%
\input\jobname.aux
\immediate\openout\@auxout=\jobname.aux
%
%

\catcode`@=11
\def\bibliographym@rk{\bgroup}


\outer\def\bye{     \par\vfill\supereject\end}

\def\og{\leavevmode\raise.3ex\hbox{$\scriptscriptstyle
\langle\!\langle\,$}}
\def \fg {\leavevmode\raise.3ex\hbox{$\scriptscriptstyle
\rangle\!\rangle\,\,$}}

\def\K{{\bf {K}}}

\def\Q{{\bf {Q}}}

\def\Z{{\bf Z}}
\def\R{{\bf R}}

\def\ovQ{{\overline \Q}}

\def\ux{{\bf x}}  \def\uz{{\bf z}}  \def\uX{{\bf X}}

\magnification=1200

\frenchspacing

\def\house#1{\setbox1=\hbox{$\,#1\,$}%
\dimen1=\ht1 \advance\dimen1 by 2pt \dimen2=\dp1 \advance\dimen2 by 2pt
\setbox1=\hbox{\vrule height\dimen1 depth\dimen2\box1\vrule}%
\setbox1=\vbox{\hrule\box1}%
\advance\dimen1 by .4pt \ht1=\dimen1
\advance\dimen2 by .4pt \dp1=\dimen2 \box1\relax}

  \def\eps{{\varepsilon}}

\def\uX{{\bf {X}}}

\def\build#1_#2^#3{\mathrel{\mathop{\kern 0pt#1}\limits_{#2}^{#3}}}

\font\fivegoth=eufm5 \font\sevengoth=eufm7 \font\tengoth=eufm10

\newfam\gothfam \scriptscriptfont\gothfam=\fivegoth
\textfont\gothfam=\tengoth \scriptfont\gothfam=\sevengoth

\def\cc{{\goth c}}

\def\smallsquare{\vbox{\hrule\hbox{\vrule height 1 ex\kern 1 ex\vrule}\hrule}}
\def\cqfd{\hfill \smallsquare\vskip 3mm}
\def\Card{{\rm Card}\,}
\def\nbdc{{\rm nbdc}}  \def\cS{{\cal S}}


\def\kdots{,\ldots ,}
\def\hfb{\hfill\break}

\def\HH{{\cal H}}
\def\LL{{\cal L}}
\def\SS{{\cal S}}

\def\E{{\bf E}}

\def\det{{\rm det}}

\def\MK{M_{{\bf K}}}
\def\ME{M_{{\bf E}}}
\def\MQ{M_{{\bf Q}}}

\def\cc{{\bf c}}
\def\xx{{\bf x}}

\def\rr{{\bf r}}

\def\HQLc{H_{Q,{\cal L},{\bf c}}}
\def\lone{\lambda_1(Q)}
\def\ltwo{\lambda_2(Q)}
\def\li{\lambda_i(Q)}

\def\indP{{\rm Ind}(P;{\bf r};{\bf x}_1\kdots {\bf x}_m)}

\def\X{{\bf X}}

\def\qed{\hfill \smallsquare\vskip 3mm}
\def\proof{\vskip2mm\noindent {\it Proof. }}


\vskip 4mm

\centerline{\bf On two notions of complexity of algebraic numbers}

\vskip 8mm
\centerline{Y{\sevenrm ANN} B{\sevenrm UGEAUD} (Strasbourg) \ {\it \&}
\ J{\sevenrm AN}-H{\sevenrm ENDRIK} E{\sevenrm VERTSE} (Leiden) \footnote{}{\rm
}}

{\narrower\narrower
\vskip 12mm

\proclaim Abstract. {
We derive new, improved lower bounds for the block complexity of an irrational
algebraic number and for the number of digit changes in the $b$-ary
expansion of an irrational algebraic number. To this end,
we apply a version of the
Quantitative Subspace Theorem by Evertse and
Schlickewei \cite{EvSc02}, Theorem 2.1.
}

}

\vskip 15mm

\centerline{\bf 1. Introduction}

\vskip 7mm

Throughout the present paper, $b$ always denotes an integer $\ge 2$
and $\xi$ is a real number with
$0 < \xi <1$. There exists a unique infinite sequence
${\bf a}=(a_j)_{j \ge 1}$
of integers from $\{0, 1, \ldots, b-1\}$, called the
$b$-ary expansion of $\xi$, such that
$$
\xi = \sum_{j\ge 1} \, {a_j \over b^j},
$$
and ${\bf a}$ does not terminate in an infinite string of the digit $b-1$.
Clearly, the sequence ${\bf a}$ is
ultimately periodic if, and only if, $\xi$
is rational. With a slight abuse
of notation, we also denote by ${\bf a}$ the infinite word
$a_1 a_2 \ldots$
To measure the complexity of $\xi$, we measure the
complexity of ${\bf a}$. Among the different ways to do this,
two notions of complexity have been recently
studied. A first one, namely the block complexity,
consists in counting the
number $p(n, \xi, b) = p(n, {\bf a})$ of distinct blocks of length
$n$ occurring in the word ${\bf a}$, that is,
$$
p(n, \xi, b) = \Card \{ a_{k+1} a_{k+2} \ldots a_{k+n} : k \ge 0 \}.
$$
A second one deals with the asymptotic behaviour
of the number of digit changes in ${\bf a}$.
The function $\nbdc$, `number of digit changes',
introduced in \cite{BuPad}, is defined by
$$
\nbdc(n, \xi, b) = \Card \{ 1 \le k \le n : a_k \not= a_{k+1} \},
\quad \hbox{for $n \ge 1$}.
$$

Suppose from now on that $\xi$ is algebraic and irrational.
Non-trivial lower bounds for $p(n, \xi, b)$ and $\nbdc(n, \xi, b)$
were obtained in \cite{AdBu07,BuPad} by means of transcendence criteria that
ultimately depend on the Schmidt Subspace Theorem \cite{SchmLN}
or on the Quantitative Roth Theorem \cite{Ro55,Lo99}. 
Respectively,
it is known that
$$
\lim_{n\to + \infty}\, {p(n, \xi, b) \over n} = +\infty  \eqno (1.1)
$$
and
$$
\nbdc(n, \xi, b) \ge 3 \, (\log n)^{1 + 1/( \omega (b) + 4)} \cdot
(\log \log n)^{-1/4},   \eqno (1.2)
$$
for every sufficiently large $n$, where $\omega(\ell)$ counts the
number of distinct prime factors of the integer $\ell$.

Both (1.1) and (1.2) are very far from what can be expected if one believes
that, regarding these notions of complexity, algebraic irrational
numbers behave like almost all real numbers (in the sense of the
Lebesgue measure). 
Thus, it is widely believed that the functions
$n \mapsto p(n, \xi, b)$ and $n \mapsto \nbdc(n, \xi, b)$
should grow, respectively, exponentially in $n$ and linearly in $n$.


The main purpose of the present paper is to
improve (1.2) for all $n$ and (1.1) for infinitely many $n$.
Our results imply that
$$
p(n, \xi , b) \ge n (\log n)^{0.09}\ \ \hbox{for infinitely many $n$}
\eqno (1.3)
$$
and
$$
\nbdc(n, \xi, b) \ge c(d) \, (\log n)^{3/2}\cdot (\log\log n)^{-1/2}
$$
for every sufficiently large $n$, where $c(d)$ is a constant depending only
on the degree $d$ of $\xi$.
In particular, we have been able to remove the dependence
on $b$ in (1.2).

The new ingredient in the proof of (1.3) is the use of a quantitative version
of the Subspace Theorem,
while (1.1)
was established by means of a standard qualitative version of the Subspace
Theorem.
Originally, quantitative versions of the Subspace Theorem 
were stated for 
a single inequality with a product of linear forms,
and then the resulting upper bound for the number of subspaces depended
on the number of places involved.
Instead, we use a version for systems of inequalities each involving
one linear form giving an upper bound for the number of subspaces
independent of the number of places.
In fact, for many applications, the version for systems
of inequalities suffices, and it leads to much better results
when many non-Archimedean places are involved.

Our paper is organized as follows. We begin by stating and discussing
our 
result
upon (1.1) in Section 2 and that upon (1.2)
in Section 3. Then, we state in Section 4 our main auxiliary tool, namely
the Quantitative Parametric Subspace Theorem from \cite{EvSc02}.
We have included an improvement of the two-dimensional case of the latter
which is needed for our improvement upon (1.2); the proof
of this improvement is included in an appendix at the end of
our paper. This Quantitative Parametric Subspace Theorem
is a statement about classes of twisted heights parametrized by
a parameter $Q$, and one can deduce from this suitable versions
of the Quantitative Subspace Theorem, dealing with (systems of)
Diophantine inequalities.
In Section 5 we deduce
a quantitative result for systems of inequalities (Theorem 5.1) 
fine-tuned for the applications in our present paper.
In the particular case where we have only two unknowns we
obtain a sharper quantitative version of a Ridout type theorem (Corollary 5.2).
The proof of Theorem 2.1 splits in Sections 6 and 7, and that of
Theorem 3.1 is given in Section 8.
Finally, further applications of our results
are discussed in Section 9.

\vskip 5mm

\centerline{\bf 2. Block complexity
of $b$-ary expansions of algebraic numbers}

\vskip 7mm

We keep the notation from the Introduction.
Recall that the real
number $\xi$ is called {\it normal in base $b$} if,
for any positive integer $n$,
each one of the $b^n$ words of length $n$
on the alphabet $\{0, 1, \ldots, b-1\}$
occurs in the $b$-ary expansion of $\xi$ with the same frequency $1/b^n$.
The first explicit example of
a number normal in base $10$, namely the number
$$
0.1234567891011121314\ldots,  
$$
whose sequence of digits is the concatenation of the sequence of all
positive integers ranged in increasing order,
was given in 1933 by Champernowne \cite{Cha33}.
It follows from the Borel--Cantelli lemma that almost all real
numbers (in the sense of the Lebesgue measure) are normal in every integer
base, but proving that a specific number, like $e$, $\pi$ or $\sqrt{2}$
is normal in some base remains a challenging open
problem. However, it is believed that
every real irrational algebraic number is normal in
every integer base. This problem, which was first formulated by
\'Emile Borel \cite{Bor50}, is likely to be very difficult.

Assume from now on that $\xi$ is algebraic and irrational.
In particular, ${\bf a}$ is not ultimately periodic.
By a result of Morse and Hedlund \cite{MoHe38,MoHe40},
every infinite word ${\bf w}$ that is not ultimately periodic
satisfies $p(n, {\bf w}) \ge n+1$ for $n \ge 1$.
Consequently, $p(n, \xi, b)\ge n+1$ holds for
every positive integer $n$. This lower bound was subsequently improved upon
in 1997 by Ferenczi and Mauduit \cite{FeMa97}, who applied a
non-Archimedean extension of Roth's Theorem
established by Ridout \cite{Rid57} to show
(see also \cite{All00}) that
$$
\lim_{n \to + \infty} \,
\bigl( p(n, \xi, b) - n \bigr) = + \infty.
$$

Then, a new combinatorial
transcendence criterion proved with the help of the
Schmidt Subspace Theorem by Adamczewski, Bugeaud, and Luca \cite{AdBuLu}
enabled Adamczewski and Bugeaud \cite{AdBu07} to establish that
$$
\lim_{n\to + \infty}\, {p(n, \xi, b) \over n} = +\infty.  \eqno (2.1)
$$

By combining ideas from \cite{BuCug} with a suitable
version of the Quantitative Subspace Theorem, we are
able to 
prove the following concerning (2.1).

\proclaim Theorem 2.1.
Let $b \ge 2$ be an integer and $\xi$ an algebraic irrational number
with $0 < \xi < 1$.
Then, for any real number $\eta$
such that $\eta < 1/11$, we have
$$
\limsup_{n\to + \infty}\, {p(n, \xi, b)
\over n (\log n)^{\eta}} = +\infty. \eqno (2.2)
$$

Ideas from \cite{BuCug} combined with Theorem 3.1 from \cite{EvSc02}
allow  
us to prove a weaker version of Theorem 2.1, namely to
conclude that (2.2) holds for any $\eta$ smaller than
$1/(4 \omega(b) + 15)$. The key point for removing
the dependence on $b$ is the use of Theorem 5.1 below,
and more precisely the fact that the exponent on $\varepsilon^{-1}$ 
in (5.9) does not depend of the cardinality of the set of places $S$. 

Remark that Theorem 2.1 does not follow from (2.1). Indeed,
there exist infinite words ${\bf w}$
having a complexity function $p$ satisfying
$$
\lim_{n\to + \infty} \, {p(n, {\bf w}) \over n} = +\infty
\quad {\rm and} \quad
\lim_{n\to + \infty}\, {p(n, {\bf w}) \over n (\log \log n)} < +\infty.
\eqno (2.3)
$$
In particular, there exist morphic words satisfying (2.3).
We refer the reader to \cite{AdBu07} for the definition of a morphic word.
An open question posed in \cite{AdBu07} asked
whether the $b$-ary expansion of an irrational algebraic number
can be a morphic word. Theorem 2.1 above allows us to make a small
step towards a negative answer. Indeed, by a result of
Pansiot \cite{Pan84}, the complexity of a
morphic word that is not ultimately periodic
is either of order $n$, $n \log \log n$, $n \log n$, or $n^2$.
It immediately follows from Theorem 2.1 that, regardless of the base $b$,
if the $b$-ary expansion of an irrational algebraic number
is generated by a morphism, then the complexity of this morphism is
either of order $n \log n$, or of order $n^2$.
However, by using 
combinatorical properties of morphic words and the 
transcendence criterion from \cite{AdBuLu}, Albert \cite{AlbTh}, 
on page 59 of his thesis, was able to show a stronger result, namely 
to prove that, regardless of the base $b$,
if the $b$-ary expansion of an irrational algebraic number
is generated by a morphism, then its complexity 
is of order $n^2$.

Note that our method yields the existence of a positive $\delta$
such that
$$
\limsup_{n\to + \infty}\, {p(n, \xi, b) \cdot (\log \log n)^{\delta}
\over n (\log n)^{1/11}} = +\infty.  \eqno (2.4)
$$
In order to avoid
painful technical details, we decided not to give a proof of (2.4).

\vskip 5mm

\centerline{\bf 3. Digit changes in $b$-ary expansions of algebraic numbers}

\vskip 7mm

Our next result is a new lower bound for the number of digit changes
in $b$-ary expansions of irrational algebraic numbers.

\proclaim Theorem 3.1.
Let $b \ge 2$ be an integer.
Let $\xi$ be an irrational, real algebraic number $\xi$
of degree $d$. There exist
an effectively computable absolute constant $c_1$
and an effectively computable constant $c_2(\xi, b)$, depending
only on $\xi$ and $b$, such that
$$
\nbdc(n, \xi, b) \ge c_1 \, {(\log n)^{3/2} \over 
(\log\log n)^{1/2} \, (\log 6d)^{1/2}}
$$
for every integer $n \ge c_2(\xi, b)$.

Theorem 3.1 improves upon Theorem 1 from \cite{BuPad},
where the exponent of $(\log n)$ depends on $b$ and
tends to $1$ as the number of prime factors of $b$ tends to
infinity.
This improvement is a consequence of the use
of the two-dimensional case of Theorem 5.1 (dealing with systems
of inequalities) 
instead 
of a result of Locher \cite{Lo99} 
(dealing with one inequality with a product of linear forms). 

Theorem 3.1 allows us to improve upon straightforwardly
many of the results from \cite{BuPad}.
We restrict our attention to Section 7 from \cite{BuPad}, that is,
to the study of the gap series 
$$
\xi_{{\bf n}, b} = \sum_{j \ge 1} \, b^{-n_j}
$$
for a given integer $b \ge 2$ and 
a non-decreasing sequence of positive integers ${\bf n} = (n_j)_{j \ge 1}$.
As mentioned in \cite{BuPad}, it 
easily follows from Ridout's Theorem \cite{Rid57} that the assumption
$$
\limsup_{j \to + \infty} \, {n_{j+1} \over n_j} > 1
$$
implies the transcendence of $\xi_{{\bf n}, b}$, see e.g. 
Satz 7 from Schneider's monograph \cite{Schn}.

In particular, for any positive real number
$\eps$, the real number $\xi_{{\bf n}, b}$
is transcendental when $n_j = 2^{[\eps j]}$, 
where $[ \, \cdot \, ]$ denotes the integer part function.
A much sharper statement, that improves Corollary 4 from \cite{BuPad}, 
follows at once from Theorem 3.1.

\proclaim Corollary 3.2.
Let $b \ge 2$ be an integer.
For any real number $\eta > 2/3$, the sum of the series
$$
\sum_{j \ge 1} \, b^{-n_j}, \quad 
\hbox{where $n_j = 2^{[j^{\eta}]}$ for $j \ge 1$,}
$$
is transcendental.

To establish Corollary 3.2, it is enough
to check that the number of positive integers $j$
such that $2^{[j^{\eta}]} \le N$ is less than some
absolute constant times $(\log N)^{1/\eta}$, and
to apply Theorem 3.1 to conclude. Stronger
transcendence results for the gap series $\xi_{{\bf n}, 2}$ follow from
\cite{BBCP,Ri06}, including the fact that Corollary 3.2 holds for any
positive $\eta$ when $b=2$.

Further results are given in Section 9.

\vskip 5mm

\centerline{\bf 4. The Quantitative Parametric Subspace Theorem}

\vskip 7mm

We fix an algebraic closure $\ovQ$ of $\Q$; all algebraic number fields
occurring henceforth will be subfields of $\ovQ$.

We introduce the necessary absolute values.
The set of places $M_{\Q}$ of $\Q$ may be identified with
$\{\infty\}\cup\{ {\rm primes}\}$. We denote by $|\cdot |_{\infty}$
the ordinary (Archimedean) absolute value on $\Q$ and for a prime $p$
we denote by $|\cdot |_p$ the $p$-adic absolute value, normalized such
that $|p|_p=p^{-1}$.

Let $\K$ be an algebraic number field. We denote by $M_{\K}$ the set of places
(equivalence classes of non-trivial absolute values) of $\K$.
The completion of $\K$ at a place $v$ is denoted by $\K_v$.
Given a place $v\in M_{\K}$, we denote by $p_v$ the place in $M_{\Q}$
lying below $v$. We choose the absolute value $|\cdot |_v$ in $v$
in such a way that the restriction of $|\cdot |_v$ to $\Q$ is $|\cdot |_{p_v}$.
Further, we define the normalized absolute value $\|\cdot\|_v$ by
$$
\|\cdot\|_v := |\cdot |_v^{d(v)}\quad\hbox{where }
d(v):={[\K_v:\Q_{p_v}]\over [\K :\Q ]}.
\eqno (4.1)
$$
These absolute values satisfy the product formula
$$ 
\prod_{v\in M_{\K}} \| x\|_v=1, \quad \hbox{for $x\in \K^*$.}
$$
Further, they satisfy the extension formula:
Suppose that $\E$ is a finite extension of $\K$ and normalized absolute values
$\|\cdot \|_w$ $(w\in \ME )$ are defined in precisely the same manner as
those for $\K$. Then if $w\in \ME$ and $v\in \MK$ is the place below $w$,
we have
$$
\| x\|_w =\| x\|_v^{d(w|v)}\ \hbox{for $x\in \K$,}
\quad\hbox{where } d(w|v):={[\E_w :\K_v]\over [\E :\K]}.
\eqno (4.2)
$$
Notice that
$$
\sum_{w|v} d(w|v)=1
\eqno (4.3)
$$
where by `$w|v$'  we indicate that $w$ runs through all places of $\E$
lying above $v$.

Let again $\K$ be an algebraic number field, and $n$ an integer $\geq 2$.
Let ${\cal L}=( L_{iv}:\, v\in M_{\K},\, i=1\kdots n)$ be a tuple of linear
forms with the following properties:
$$
\eqalignno{
&L_{iv}\in \K [X_1\kdots X_n]\ \hbox{for $v\in\MK$, $i=1\kdots n$,}&(4.4)
\cr
&L_{1v}=X_1\kdots L_{nv}=X_n\ \hbox{for all but finitely many $v\in M_{\K}$,}
&(4.5)
\cr
&\det (L_{1v}\kdots L_{nv})=1\ \hbox{for } v\in M_{\K},&(4.6)
\cr
&\Card\Big(\bigcup_{v\in M_{\K}}\{ L_{1v}\kdots L_{nv}\}\Big)\leq r.&(4.7)
}
$$
Further, we define
$$
\HH =\HH (\LL ) =\prod_{v\in M_{\K}}\max_{1\leq i_1<\cdots <i_n\leq s}\|
\det (L_{i_1}\kdots L_{i_n})\|_v
\eqno (4.8)
$$
where we have written $\{ L_1\kdots L_s\}$ for 
$\bigcup_{v\in M_{\K}}\{ L_{1v}\kdots L_{nv}\}$.  

Let ${\bf c}=(c_{iv}:\, v\in M_{\K},\, i=1\kdots n)$ be a tuple of reals
with the following properties:
$$
\eqalignno{
&c_{1v}=\cdots =c_{nv}=0\ \hbox{for all but finitely many $v\in M_{\K}$,}&(4.9)
\cr
&\sum_{v\in\MK}\sum_{i=1}^n c_{iv}=0,&(4.10)
\cr
&\sum_{v\in\MK}\max (c_{1v}\kdots c_{nv})\leq 1.&(4.11)
}
$$

Finally, for any finite extension $\E$ of $\K$
and any place $w\in\ME$ we define
$$
L_{iw}=L_{iv},\ c_{iw}=d(w|v)c_{iv}\quad\hbox{for } i=1\kdots n,
\eqno (4.12)
$$
where $v$ is the place of $\MK$ lying below $w$ and $d(w|v)$ is given by (4.2).

We define a so-called twisted height $H_{Q,\LL ,\cc}$ on $\ovQ^n$ as follows.
For $\xx\in\K^n$ define
$$
\HQLc (\xx ):=\prod_{v\in\MK}\max_{1\leq i\leq n}\| L_{iv}(\xx )\|_vQ^{-c_{iv}}.
$$
More generally, for $\xx\in\ovQ^n$ take any finite extension $\E$ of $\K$
with $\xx\in\E^n$ and put
$$
\HQLc (\xx ):=\prod_{w\in\ME}\max_{1\leq i\leq n}\| L_{iw}(\xx )\|_wQ^{-c_{iw}}.
\eqno (4.13)
$$
Using (4.12), (4.2), (4.3), and basic properties of degrees of field
extensions, one easily shows that this does not depend on the choice of $\E$.

\proclaim Proposition 4.1.
Let $n$ be an integer $\geq 2$, let $\LL =(L_{iv}:\, v\in\MK ,\, i=1\kdots n)$
be a tuple of linear forms
satisfying (4.4)--(4.7) and $\cc =(c_{iv}:\, v\in\MK ,\, i=1\kdots n)$
a tuple of reals satisfying (4.9)--(4.11). Further, let $0<\delta \leq 1$.
\hfb
Then there are proper linear subspaces $T_1\kdots T_{t_1}$ of $\ovQ^n$,
all defined over $\K$, with
$$
t_1=t_1(n,r,\delta )=
\left\{
\matrix{
4^{(n+8)^2}\delta^{-n-4}\log (2r)\log\log(2r)\hfill &\hbox{if } n\geq 3,\hfill
\cr
2^{25}\delta^{-3}\log (2r)\log\big(\delta^{-1}\log (2r)\big)\hfill
&\hbox{if } n=2\hfill
}
\right.
$$
such that the following holds: for every real $Q$ with
$$
Q>\max \Big(\HH^{1/{r\choose n}},n^{2/\delta}\Big)
$$
there is a subspace $T_i\in\{ T_1\kdots T_{t_1}\}$ such that
$$
\{ \xx\in\ovQ^n :\, \HQLc (\xx )\leq Q^{-\delta}\}\subset T_i\, .
$$

For $n\geq 3$ this is precisely Theorem 2.1 of \cite{EvSc02},
while for $n=2$ this is an improvement of this theorem.
This improvement can be obtained by combining some lemmata from \cite{EvSc02}
with
more precise computations in the case $n=2$.
We give more details in the appendix at the end of the present paper.

\vskip 5mm

\centerline{{\bf 5. Systems of inequalities}}

\vskip 7mm

For every place $p\in\MQ =\{\infty\}\cup\{{\rm primes}\}$ 
we choose an extension
of $|\cdot |_p$ to $\ovQ$ which we denote also by $|\cdot |_p$.
For a linear form $L=\sum_{i=1}^n \alpha_iX_i$ with coefficients in $\ovQ$
we define the following:
We denote by $\Q (L)$ the field generated by the coefficients of $L$, i.e.,
$\Q (L):= \Q (\alpha_1\kdots\alpha_n)$; for any map $\sigma$ from $\Q (L)$ to
any other field we define $\sigma (L):=\sum_{i=1}^n \sigma (\alpha_i)X_i$;
and the inhomogeneous height of $L$ is given by
$H^*(L):=\prod_{v\in\MK} \max (1,\|\alpha_1\|_v\kdots \|\alpha_n\|_v)$,
where $\K$ is any number field containing $\Q (L)$.
Further, we put $\| L\|_v:=\max(\|\alpha_1\|_v\kdots \|\alpha_n\|_v)$
for $v\in\MK$.

Let $n$ be an integer with $n\geq 2$, $\varepsilon$ a real with
$\varepsilon >0$ and $S=\{ \infty ,p_1\kdots p_t\}$ a finite subset of
$\MQ$ containing the infinite place. Further, let $L_{ip}$ $(p\in S,\,
i=1\kdots n)$ be linear forms in $X_1\kdots X_n$ with coefficients in $\ovQ$
such that
$$
\eqalignno{
&\det (L_{1p}\kdots L_{np})=1\ \hbox{for $p\in S$,}
&(5.1)
\cr
&\Card\Big(\bigcup_{p\in S}\{ L_{1p}\kdots L_{np}\}\Big)\leq R,&(5.2)
\cr
&[\Q (L_{ip}):\Q ]\leq D\ \hbox{for $p\in S,\, i=1\kdots n$,}&(5.3)
\cr
&H^*(L_{ip})\leq H\ \hbox{for $p\in S,\, i=1\kdots n$,}&(5.4)
\cr
}
$$
and $e_{ip}$ $(p\in S,\,i=1\kdots n)$ be reals satisfying
$$
\eqalignno{
&e_{i\infty}\leq 1\ (i=1\kdots n),\ \
e_{ip}\leq 0\ (p\in S\setminus\{\infty\},\, i=1\kdots n),
&(5.5)
\cr
&\sum_{p\in S}\sum_{i=1}^n e_{ip} = -\varepsilon .&(5.6)
}
$$
Finally let $\Psi$ be a function from $\Z^n$ to $\R_{\geq 0}$. We consider
the system of inequalities
$$
\eqalign{
|L_{ip}(\xx )|_p\leq \Psi (\xx )^{e_{ip}}\ &(p\in S,\, i=1\kdots n)\cr
&
\hbox{in $\xx\in\Z^n$ with $\Psi (\xx )\not= 0$.}
}
\eqno (5.7)
$$

\proclaim Theorem 5.1.
The set of solutions of (5.7) with
$$
\Psi (\xx )>\max \big( 2H, n^{2n/\varepsilon}\big)
\eqno (5.8)
$$
is contained in the union of at most
$$
\left\{
\matrix{
8^{(n+9)^2}(1+\varepsilon^{-1})^{n+4}\log (2RD)\log\log (2RD)\hfill&
\hbox{if $n\geq 3$}\hfill\cr
2^{32}(1+\varepsilon^{-1})^3
\log (2RD)\log\big( (1+\varepsilon^{-1})\log (2RD)\big)\hfill&
\hbox{if $n=2$}\hfill
}
\right.
\eqno (5.9)
$$
proper linear subspaces of $\Q^n$.

\vskip 3mm  
\noindent
{\it Remark.}
Let $\|\cdot\|$ be any vector norm on $\Z^n$. 
Then for the solutions $\xx$ of (5.7) we have, 
in view of (5.5),
$$
\|\xx\| \ll \max_{1\leq i\leq n} |L_{i\infty}(\xx )|\ll \Psi (\xx ).
$$
So it would not have been a substantial restriction
if in the formulation of Theorem 5.1 
we had restricted the function $\Psi$ to vector norms.
But for applications it is convenient to allow other functions for $\Psi$.

\vskip2mm

We deduce from Theorem 5.1 a quantitative Ridout type theorem.
Let $S_1$, $S_2$ be finite, possibly empty sets of prime numbers,
put $S:=\{\infty\}\cup S_1\cup S_2$,
let $\xi\in\ovQ$ be an algebraic number,
let $\varepsilon >0$,
and let $f_p$ $(p\in S)$ be reals such that
$$
f_p\geq 0\ \hbox{for $p\in S$,}\ \ \sum_{p\in S} f_p=2+\varepsilon .
\eqno (5.10)
$$

We consider the system of inequalities
$$
\left\{
\matrix{
|\xi -{x\over y}|\leq y^{-f_{\infty}},\hfill\cr
|x|_p\leq y^{-f_p}\ (p\in S_1)\hfill\cr
|y|_p\leq y^{-f_p}\ (p\in S_2)\hfill
}
\right\}
\hbox{in $(x,y)\in\Z^2$ with $y>0$.}
\eqno (5.11)
$$
Define the height of $\xi$ by
$H(\xi ):=\prod_{v\in\MK}\max (1,\|\xi \|_v)$ where $\K$ is any algebraic
number field with $\xi\in\K$. Suppose that $\xi$ has degree $d$.

\proclaim Corollary 5.2.
The set of solutions of (5.11) with
$$
y> \max\big( 2H(\xi ), 2^{4/\varepsilon}\Big)
\eqno (5.12)
$$
is contained in the union of at most
$$
2^{32}(1+\varepsilon^{-1})^3
\log (6d)\log\big( (1+\varepsilon^{-1})\log (6d)\big)
\eqno (5.13)
$$
one-dimensional linear subspaces of $\Q^2$.

To obtain Corollary 5.2 one simply has to apply Theorem 5.1 with $n=2$,
$S=\{\infty\}\cup S_1\cup S_2$ and with
$$
\eqalign{
&L_{1\infty}=X_1-\xi X_2,\, L_{2\infty}=X_2,\cr
&L_{1p}=X_1,\, L_{2p}=X_2\ \hbox{for $p\in S_1\cup S_2$,}\cr
&e_{1\infty}=1-f_{\infty},\, e_{2\infty}=1,\cr
&e_{1p}=-f_p,\, e_{2p}=0\ \hbox{for $p\in S_1$,}\cr
&e_{1p}=0,\, e_{2p}=-f_p\ \hbox{for $p\in S_2$,}\cr
&\Psi (\xx )=|x_2|\ \hbox{for $\xx =(x_1,x_2)\in\Z^2$.}
}
$$
It is straightforward to verify that (5.1) is satisfied,
and that (5.2), (5.3), (5.4) are satisfied with $R=3$, $D=d$,
$H=H(\xi )$, respectively.
Further, it follows at once from (5.10) that
(5.5) and (5.6) are satisfied.
\vskip 5mm

\noindent
{\it Proof of Theorem 5.1.}
Let $\K$ be a finite normal extension of $\Q$, containing the coefficients
of $L_{ip}$ as well as the conjugates over $\Q$ of these coefficients,
for $p\in S$, $i=1\kdots n$. For $v\in\MK$ we put $d(v):=[\K_v:\Q_{p_v}]$
where $p_v$ is the place of $\Q$ below $v$, and
$$
s(v)=d(v)\ \hbox{if $v$ is Archimedean},\ \
s(v)=0\ \hbox{if $v$ is non-Archimedean}.
$$
Recall that every $|\cdot |_p$ $(p\in\MQ )$ has been extended to $\ovQ$
so in particular to $\K$.
For every $v\in\MK$ there is an automorphism $\sigma_v$ of $\K$
such that $|\sigma_v(\cdot )|_p$ represents $v$. So by (4.1) we have
$$
\| x\|_v=|\sigma_v(x)|_{p_v}^{d(v)}\ \hbox{for $x\in\K$.}
\eqno (5.14)
$$
Let $T$ denote the set of places of $\K$ lying above the places in $S$.
Define linear forms $L_{iv}$ and reals $e_{iv}$ ($v\in\MK,\, i=1\kdots n$) by
$$
L_{iv}=\sigma_v^{-1}(L_{i,p_v})\ (v\in T),\ \
L_{iv}=X_i\ (v\in \MK\setminus T)\eqno (5.15)
$$
and
$$
e_{iv}=d(v)e_{i,p_v}\ (v\in T),\ \ e_{iv}=0\ (v\in \MK\setminus T),\eqno (5.16)
$$
respectively. Then system (5.7) can be rewritten as
$$
\eqalign{
\| L_{iv}(\xx )\|_v\leq \Psi (\xx )^{e_{iv}}\ &(v\in\MK,\, i=1\kdots n)\cr
&\hbox{in $\xx\in\Z^n$ with $\Psi (\xx )\not= 0$.}
}
\eqno (5.17)
$$
Notice that in view of (5.17), (5.5), (5.6), and
$\sum_{v|p} d(v)=1$ for $p\in\MQ$
we have
$$
e_{iv}\leq s(v)\ (i=1\kdots n),\ \
\sum_{v\in \MK}\sum_{i=1}^n e_{iv}\leq -\varepsilon .
\eqno (5.18)
$$
Further, by (5.2), (5.15),
$$
\Card\bigcup_{v\in\MK}\{ L_{1v}\kdots L_{nv}\} \leq r:= n+DR.
\eqno (5.19)
$$

Now define
$$
\delta := {\varepsilon\over n+\varepsilon},\eqno (5.20)
$$
let $\LL =(L_{iv}:\, v\in\MK ,\, i=1\kdots n)$, and define the tuple
of reals $\cc =(c_{iv}:\, v\in\MK ,\, i=1\kdots n)$ by
$$
c_{iv}:= \big( 1+(\varepsilon /n)\big)^{-1}\Big( e_{iv}-{1\over n}\sum_{j=1}^n e_{jv}\Big).
\eqno (5.21)
$$
Let $\HH =\HH (\LL )$ be the quantity defined by (4.8) and $\HQLc$ the twisted
height defined by (4.13).
We want to apply Proposition 4.1, and to this end we have to verify
the conditions (4.4)--(4.7) and (4.9)--(4.11).
Condition (4.4) is obvious.
(5.1) and (5.15) imply (4.5),(4.6), while (4.7)
is (5.18).
Condition (4.9) is satisfied in view of (5.16), (5.20), while (4.10)
follows at once from (5.21). To verify (4.11), observe that by
(5.21), (5.18) we have
$$
\eqalign{
&\sum_{v\in\MK}\max (c_{1v}\kdots c_{nv})\cr
&\qquad
\leq\Big( 1+{\varepsilon\over n}\Big)^{-1}\Big( \sum_{v\in\MK} s(v) -
{1\over n}\sum_{v\in\MK}\sum_{j=1}^n e_{jv}\Big)\cr
&\qquad
=\Big( 1+{\varepsilon\over n}\Big)^{-1}\Big( 1+{\varepsilon\over n}\Big)=1.
}
$$

The following lemma connects system (5.7) to Proposition 4.1.

\proclaim Lemma 5.3.
Let $\xx$ be a solution of (5.7) with (5.8). Put
$$
Q := \Psi (\xx )^{1+\varepsilon /n}.
$$
Then
$$
\eqalignno{
&\HQLc (\xx )\leq Q^{-\delta},&(5.22)
\cr
&Q\geq\max \Big(\HH^{1/{r\choose n}},n^{2/\delta}\Big).&(5.23)
}
$$

\noindent{\it Proof.}
As observed above, $\xx$ satisfies (5.17). In combination
with (5.21) this yields
$$
\eqalign{
\| L_{iv}(\xx )\|_vQ^{-c_{iv}}&=
\| L_{iv}(\xx )\|_v\cdot \Psi (\xx )^{-e_{iv}}
\cdot\Psi (\xx )^{{1\over n}\sum_{j=1}^n e_{jv}}\cr
&\leq\Psi (\xx )^{{1\over n}\sum_{j=1}^n e_{jv}}
}
$$
for $v\in\MK$, $i=1\kdots n$. By taking the product over $v\in\MK$
and using (5.18), (5.20) we obtain
$$
\eqalign{
\HQLc (\xx) &=
\prod_{v\in\MK}\max_{1\leq i\leq n} \| L_{iv}(\xx )\|_vQ^{-c_{iv}}
\cr
&\leq\Psi (\xx)^{-\varepsilon /n}=Q^{-\delta}.
}
$$
This proves (5.22).

To prove (5.23), write
$\bigcup_{v\in\MK}\{ L_{1v}\kdots L_{nv}\}=\{ L_1\kdots L_s\}$.
Then $s\leq r$ by (5.19).
By (5.4), (5.15) we have $H^*(L_{iv})\leq H$ for $v\in \MK$, $i=1\kdots n$.
By applying e.g., Hadamard's inequality for the Archimedean places and
the ultrametric inequality for the non-Archimedean places,
we obtain for $i_1\kdots i_n\in\{ 1\kdots s\}$, $v\in\MK$, 
$$
\eqalign{
\|\det (L_{i_1}\kdots L_{i_n})\|_v &\leq (n^{n/2})^{s(v)}
\prod_{j=1}^n\| L_{i_j}\|_v\cr
&\leq
(n^{n/2})^{s(v)}
\prod_{i=1}^s\max (1,\| L_i\|_v),\cr
}
$$
hence
$$
\HH \leq n^{n/2}\prod_{i=1}^r H^*(L_i)\leq n^{n/2}H^r.
$$
Together with (5.19), (5.20) this implies
$$
\eqalign{
\max\Big(\HH^{1/{r\choose n}},n^{2/\delta}\Big)&\leq
\max\Big( n^{n/2{r\choose n}}H^{r/{r\choose n}}, n^{2(n+\varepsilon )/\varepsilon}\Big)\cr
&\leq \max\Big( 2H ,n^{2n/\varepsilon}\Big)^{1+\varepsilon /n}.
}
$$
So if $\xx$ satisfies (5.8), then $Q=\Psi (\xx )^{1+\varepsilon /n}$
satisfies (5.23). This proves Lemma 5.3.\cqfd
\vskip5mm

We apply Proposition 4.1 with the values
of $r,\delta$ given by (5.19), (5.20),
i.e., $r=n+DR$ and $\delta ={\varepsilon \over n+\varepsilon}$.
It is straightforward to show that for these choices of $r,\delta$
the quantity $t_1$ from Proposition 4.1
is bounded above by the quantity in (5.9).
By Proposition 4.1, there are proper linear subspaces $T_1\kdots T_{t_1}$
of $\ovQ^n$ such that for every $Q$ with (5.23) there is
$T_i\in\{ T_1\kdots T_{t_1}\}$ with
$$
\{\xx\in\ovQ^n:\, \HQLc (\xx )\leq Q^{-\delta}\}\subset T_i.
$$
Now Lemma 5.3 implies that the solutions $\xx$ of (5.7) with (5.8) lie in
$\cup_{i=1}^{t_1} (T_i\cap \Q^n)$. Theorem 5.1 follows.\cqfd

\vskip 5mm

\centerline{\bf 6. A combinatorial lemma for the proof of Theorem 2.1}

\vskip 7mm

In this section, we establish the following lemma.

\proclaim Lemma 6.1.
Let $b \ge 2$ be an integer.
Let $c$ and $u$ be positive real numbers.
Let $\xi$ be an irrational real number such that $0< \xi < 1$ and
$$
p(n, \xi, b) \le c n (\log n)^u, \quad
\hbox{for $n \ge 1$}.
$$
Then
for every
positive real number $v < u$, there exist integer sequences
$(r_n)_{n \ge 1}$, $(t_n)_{n \ge 1}$, $(p_n)_{n \ge 1}$ and
a positive real number $C$ depending only on $c,u,v$ such that
$$
|b^{t_n} \xi - b^{r_n} \xi - p_n|
\le \bigl(b^{t_n} \bigr)^{- (\log t_n)^{-v}},
$$
$$
0\leq r_n<t_n,\quad 2 t_n < t_{n + 1}, \quad
t_n \le (2n)^{C n}, \quad \hbox{for $n \ge 1$}.  \eqno (6.1)
$$
Furthermore, $b$ does
not divide $p_n$ if $r_n \ge 1$.

\noindent {\it Proof.}
Let $b$ and $\xi$ be as in the statement of the lemma.
Let ${\bf a}$ denote the $b$-ary expansion of $\xi$.
Throughout this proof, $c_1, c_2, \ldots $
are positive constants depending only on $c,u,v$.
The length of a finite word
$W$, that is, the number of letters
composing $W$, is denoted by $\vert W \vert$.
The infinite word $W^{\infty}$ is obtained by concatenation of
infinitely many copies of the finite word $W$.

By assumption, the complexity function of
${\bf a}$ satisfies
$$
p(n, {\bf a}) \le c_1 \,  n (\log n)^u, \quad \hbox{for $n \ge 1$}.
$$
Our aim is to show that
there exists a
(in some sense) `dense' sequence of rational approximations to $\xi$
with special properties.

Let $\ell \ge 2$ be an integer, and denote by $A(\ell)$
the prefix of ${\bf a}$ of length $\ell$.
By the {\it Schubfachprinzip}, there exist
(possibly empty) words $U_{\ell},
V_{\ell}, W_{\ell}$ and $X_{\ell}$ such that
$$
A(\ell) = U_{\ell} V_{\ell} W_{\ell} V_{\ell} X_{\ell},
$$
and
$$
|V_{\ell}| \ge c_2 \, \ell (\log \ell)^{-u}.
$$
Set $r_{\ell} = |U_{\ell}|$ and $s_{\ell} = |V_{\ell} W_{\ell}|$.
We choose the words $U_\ell$, $V_\ell$, $W_\ell$ and $X_\ell$
in such a way that $|V_{\ell}|$
is maximal and, among the corresponding factorisations of $A(\ell)$,
such that $|U_{\ell}|$ is minimal. In particular,
either $U_{\ell}$ is the empty word,
or the last digits of $U_{\ell}$ and $V_{\ell} W_{\ell}$ are different.

If $\xi_{\ell}$ denotes the rational number with
$b$-ary expansion $U_{\ell}
(V_{\ell} W_{\ell})^{\infty}$, then there exists an integer $p_{\ell}$
such that
$$
\xi_{\ell} = {p_{\ell} \over b^{r_\ell} (b^{s_\ell} - 1) },
\qquad
|\xi - \xi_{\ell}| \le b^{-r_{\ell} - s_{\ell} - |V_{\ell}|},
$$
and $b$ does not divide $p_{\ell}$ if $r_{\ell} \ge 1$.

Take $t_{\ell} = r_{\ell} + s_{\ell}$. Then
$$
\ell\ge t_{\ell}\ge s_{\ell}\geq c_2\ell(\log\ell )^{-u}.
\eqno (6.2)
$$
Hence,
$$
\eqalign{
|b^{t_\ell} \xi - b^{r_\ell} \xi - p_\ell|
&\le b^{-c_2\ell(\log\ell )^{-u}}
\cr
&\le (b^{t_{\ell}})^{-c_3(\log t_{\ell})^{-u}}.
}
$$

We construct a sequence of positive integers $(\ell_k)_{k=1}^{\infty}$
such that for every $k\geq 1$,
$$
\eqalignno{
&|b^{t_{\ell_k}} \xi - b^{{r_{\ell_k}}} \xi - p_{\ell_k}|
\leq (b^{t_{\ell_k}})^{-(\log t_{\ell_k})^{-v}},
&(6.3)
\cr
&t_{\ell_{k+1}}>2t_{\ell_k},&(6.4)
\cr
&t_{\ell_k}\le (2k)^{Ck}.&(6.5)
}
$$
Then a slight change of notation establishes the lemma.

Let
$\ell_1$ be the smallest positive integer $\ell$ such that
$c_3(\log t_{\ell})^u\ge (\log t_{\ell})^v$. Further, for $k=1,2,\ldots$,
let $\ell_{k+1}$ be the smallest positive integer $\ell$ such that
$t_{\ell}>2t_{\ell_k}$. This sequence is well-defined by (6.2).
It is clear that (6.3), (6.4) are satisfied.

To prove (6.5),
observe that if $\ell$ is any integer with $c_2\ell (\log\ell )^{-u}>2\ell_k$
then, by (6.2), $t_{\ell}>2\ell_k\geq 2t_{\ell_k}$. This shows that there
is a constant $c_4$ such that $\ell_{k+1}\le c_4\ell_k(\log\ell_k)^u$.
Now an easy induction yields that there exists a constant $C$,
depending only on $c,u,v$, such that $\ell_k\leq (2k)^{Ck}$ for $k\geq 1$.
Invoking again (6.2) we obtain (6.5).
\qed

\vskip 5mm


\centerline{\bf 7. Completion of the proof of Theorem 2.1}

\vskip 7mm

Let $\xi$ be an algebraic irrational real number.
Let $v$ be a real number such that $0 < v < 1/11$.
Define the positive real number $\eta$ by
$$
(11 + 2 \eta) (v + \eta) + \eta = 1. \eqno (7.1)
$$

We assume that there exists a positive constant $c$
such that the complexity function of $\xi$ in base $b$ satisfies
$$
p(n, \xi, b) \le c n (\log n)^{v + \eta}\quad
\hbox{for $n \ge 1$}, \eqno (7.2)
$$
and we will derive a contradiction. Then Theorem 2.1 follows.

Let $N$ be a sufficiently large integer.
We will often use the fact that $N$
is large, in order to absorb numerical constants.

Let $(r_n)_{n \ge 1}$,
$(t_n)_{n \ge 1}$, and $(p_n)_{n \ge 1}$ be the sequences
given by Lemma 6.1 applied with $u:= v + \eta$.
Set
$$
\eps = (\log t_N)^{-v},  \eqno (7.3)
$$
and observe that, in view of (6.1) and (7.3), we have
$$
\eps^{-1} = (\log t_N)^v \le N^{v+\eta}.  \eqno (7.4)
$$
For $n=1, \ldots , N$, we have
$$
|b^{t_n} \xi - b^{r_n} \xi - p_n| <
(b^{t_n})^{-\eps}.  \eqno (7.5)
$$
Put
$$
k := [2/\eps] + 1. \eqno (7.6)
$$
For each $n = 1, \ldots , N$ there is $\ell \in \{0, 1, \ldots , k-1\}$
such that
$$
{\ell \over k} \le {r_n \over t_n}
< {\ell + 1 \over k}.
$$
For the moment, we consider those $n \in \{1, \ldots , N\}$ such that
$$
{N \over 2} \le n \le N, \quad {\ell \over k} \le {r_n \over t_n}
< {\ell + 1 \over k}, \eqno (7.7)
$$
where $\ell \in \{0, 1, \ldots , k-1\}$ is fixed, and show that the vectors
$$
\ux_n := (b^{t_n}, b^{r_n}, p_n)
$$
satisfy a system of inequalities to which Theorem 5.1 is applicable.

Let $S = \{\infty\} \cup \{p : p \mid b\}$
be the set of places on $\Q$ composed of the infinite place
and the finite places corresponding to the prime divisors of $b$.
We choose
$$
\Psi (\xx )=x_1\ \ \hbox{for $\xx =(x_1,x_2,x_3)\in\Z^3$.}
$$
We introduce the linear forms with real algebraic coefficients
$$
L_{1 \infty} (\uX) = X_1, \quad
L_{2 \infty} (\uX) = X_2, \quad
L_{3 \infty} (\uX) = -\xi X_1 +\xi X_2 + X_3,
$$
and, for every prime divisor $p$ of $b$, we set
$$
L_{1 p} (\uX) = X_1, \quad
L_{2 p} (\uX) = X_2, \quad
L_{3 p} (\uX) = X_3.
$$
Set also
$$
e_{1 \infty} = 1, \quad e_{2 \infty} = {\ell + 1 \over k}, \quad
e_{3 \infty} = - \eps,
$$
and, for every prime divisor $p$ of $b$,
$$
e_{1 p} = {\log |b|_p \over \log p} ,
\quad e_{2 p} = {\log |b|_p \over \log p} \cdot {\ell \over k}, \quad
e_{3 p} = 0.
$$
Notice that
$$
\eqalign{
&\sum_{p \in S} \, \sum_{i=1}^3 e_{ip} =  - (\eps - 1/k), \cr
&e_{i\infty}\leq 1\ \ (i=1,2,3),\cr
&e_{ip}\leq 0\ \ (p\in S\setminus\{\infty\},\ i=1,2,3).}
\eqno (7.8)
$$
Furthermore,
$$
\det (L_{1p}, L_{2p}, L_{3p}) = 1,
\quad \hbox{for $p \in S$}. \eqno (7.9)
$$
Writing $d := [\Q(\xi) : \Q]$, we have
$$
\eqalign{
&\Card \bigcup_{p \in S} \, \{L_{1p}, L_{2p}, L_{3p}\} = 4,\cr
&[\Q(L_{ip}) : \Q] \le d, \ \
\hbox{for $p \in S$, $i = 1, 2, 3$.}
}
\eqno (7.10)
$$
Further,
$$
\max_{p \in S, i=1,2,3}^{} \, H^*(L_{ip}) = H(\xi). \eqno (7.11)
$$
(7.8)--(7.11) imply that the linear forms $L_{ip}$ and reals $e_{ip}$
defined above satisfy the conditions (5.1)--(5.6) of Theorem 5.1
with $n=3$, $R=4$, $D=d$, $H=H(\xi )$.

It is clear from (7.5),
(7.8) that for any integer $n$ with (7.7) we have
$$
|L_{ip} (\ux_n)|_p \le \Psi (\ux_n )^{e_{ip}},
\quad \hbox{for $p \in S$, $i = 1, 2, 3$.}
$$
Assuming that $N$ is sufficiently large, we infer from (6.1), (7.4) that for
every $n$ with (7.7) we have
$$
\Psi (\ux_n ) = b^{t_n} \ge 2^{2^{N/2} - 1} >
\max\{2 H(\xi), 3^{6/(\eps - 1/k)} \}.
$$
Now, Theorem 5.1 implies that
the set of vectors $\ux_n = (b^{t_n}, b^{r_n}, p_n)$ with
$n$ satisfying (7.7) is contained in the union of at most
$$
A_1 := 8^{144} \Big(1+ (\eps - 1/k)^{-1}\Big)^{-7} \log (8d)
\, \log \log (8d)
$$
proper linear subspaces of $\Q^3$.
We now consider the vectors $\ux_n$ with $N/2 \le n \le N$
and drop the condition $\ell/k \le r_n/t_n < (\ell + 1)/k$.
Then by (7.6), 
for any sufficiently large $N$, the set
of vectors $\ux_n = (b^{t_n}, b^{r_n}, p_n)$ with
$$
{N \over 2} \le n \le N,
$$
lies in the union of at most
$$
k A_1 \le (\eps^{-1})^{8 + \eta}
$$
proper linear subspaces of $\Q^3$.

We claim that if $N$ is sufficiently large,
then any two-dimensional linear subspace of $\Q^3$
contains at most
$(\eps^{-1})^{3 + \eta}$ vectors $\ux_n$. Having achieved this, it follows
by (7.1), (7.4) that
$$
{N \over 2} \le (\eps^{-1})^{8 + \eta} \cdot
(\eps^{-1})^{3 + \eta} \le N^{(11 + 2 \eta) (v + \eta)}
= N^{1 - \eta},
$$
which is clearly impossible if $N$ is sufficiently large. Thus (7.2) leads to
a contradiction.

So let $T$ be a two-dimensional linear subspace of $\Q^3$,
say given by an equation $z_1 X_1 + z_2 X_2 + z_3 X_3 = 0$
where we may assume that $z_1$, $z_2$, $z_3$ 
are integers without a common prime divisor.  
Let
$$
{\cal N} = \{i_1 < i_2 < \ldots < i_r\}
$$
be the set of $n$ with $N/2 \le n \le N$ and $\ux_n \in T$.
So we have to prove that $r \le (\eps^{-1})^{3 + \eta}$.

Recall that by Lemma 6.1, for every $n\ge 1$ we have either $r_n=0$,
or $r_n>0$ and $b$ does not divide $p_n$.
Hence the vectors $\ux_n$, $n \ge 1$, are pairwise
non-collinear. So the exterior product of $\ux_{i_1}, \ux_{i_2}$
must be a non-zero multiple of
$\uz = (z_1, z_2, z_3)$, and therefore
$$
\max\{|z_1|, |z_2|, |z_3|\} \le 2 b^{2t_{i_2}}.  \eqno (7.12)
$$
By combining (7.5) with
$z_1 b^{t_n} + z_2 b^{r_n} + z_3 p_n = 0$, eliminating $b^{r_n}$,
it follows that for $n$ in ${\cal N}$
$$
\biggl| {\xi (z_1 + z_2) \over \xi z_3 - z_2} - {-p_n \over b^{t_n}}
\biggr| < \biggl| {z_2 \over \xi z_3 - z_2} \biggr| \cdot
(b^{t_n})^{-1 -\eps}. \eqno (7.13)
$$
We want to apply Corollary 5.2 with
$\xi (z_1 + z_2) /(\xi z_3 - z_2)$ instead of $\xi$.

Recall that $d$ denotes the degree of $\xi$.
By (7.12), assuming that $N$ is sufficiently large, we have
$$
H\biggl( {\xi (z_1 + z_2) \over \xi z_3 - z_2} \biggr)
\le 4 b^{2 t_{i_2}} H(\xi) \le b^{3 t_{i_2}}.
$$
Likewise,
$$
\biggl| {z_2 \over \xi z_3 - z_2} \biggr|
\le H\biggl( {\xi (z_1 + z_2) \over \xi z_3 - z_2} \biggr)^d
\le b^{3 d t_{i_2}}.
$$
There is no loss of generality to assume that there is an integer $k\le r$
with
$$
b^{t_{i_k}} \ge b^{(3 d t_{i_2})^{2/\eps}}. \eqno (7.14)
$$
Indeed, if there is no such $k$
then we infer from (6.1) that
$$
b^{t_{i_2}^{2^{r-3}}} \le b^{(3 d t_{i_2})^{2/\eps}} 
\le b^{t_{i_2}^{4/\eps}},
$$
hence
$$
r \le 3 + {\log (4/\eps) \over \log 2},
$$
which is stronger than what we have to prove.
Letting $k_0$ be the smallest integer $k$ with (7.14), we have
$$
b^{t_{i_{k_0}}} \ge b^{(3 d t_{i_2})^{2/\eps}}, \quad
k_0 \le 4 + {\log (4/\eps) \over \log 2}.  \eqno (7.15)
$$

Let ${\cal N'} = \{i_{k_0}, i_{k_0 + 1}, \ldots , i_r\}$.
We divide this set further into
$$
{\cal N''}=\{ n\in {\cal N'}:\, r_n\not=0\},\ \
{\cal N'''}=\{ n\in {\cal N'}:\, r_n=0\}.
$$

By (7.13) we have for $n$ in ${\cal N''}$
$$
\biggl| {\xi (z_1 + z_2) \over \xi z_3 - z_2} - {-p_n \over b^{t_n}}
\biggr| < (b^{t_n})^{-1 - \eps/2}.  \eqno (7.16)
$$
Let $S_1 = \emptyset$ and $S_2 = \{p : p \mid b\}$.
Then for $\ell \in S_2$ we have
$$
|b^{t_n}|_{\ell} \le (b^{t_n})^{\log |b|_{\ell} / (\log b)}. \eqno (7.17)
$$
Lastly,
$$
b^{t_n} \ge b^{(3 d t_{i_2})^{2/\eps}}  \ge
\max\biggl\{ H\biggl( {\xi (z_1 + z_2) \over \xi z_3 - z_2} \biggr),
2^{4/\eps} \biggr\}.  \eqno (7.18)
$$
Now, (7.16), (7.17) and (7.18) imply that all the conditions of Corollary 5.2
are satisfied with $\eps /2$ instead of $\eps$ and with
$$
x=p_n,\ y=b^{t_n}, \quad
f_{\infty} = 1 + {\eps \over 2}, \quad
f_{\ell} = - {\log |b|_{\ell} \over \log b} \enspace (\ell \in S_2).
$$
Notice that
$$
f_{\infty} + \sum_{\ell \in S_2} \, f_{\ell} = 2 + \eps / 2,
$$
and $f_{\infty} \ge 0$, $f_{\ell} \ge 0$ for $\ell \in S_2$.
Consequently, the set of vectors $(p_n, b^{t_n})$, $ n \in {\cal N''}$,
lies in the union of at most
$$
B(d,\eps ):= 2^{32}\big( 1+2\eps^{-1}\big)^3\log (6d)
\log\big( (1+2\eps^{-1})\log (6d)\big)
\eqno (7.19)
$$
one-dimensional linear subspaces of $\Q^2$. But
the vectors $(p_n, b^{t_n})$, $ n \in {\cal N''}$,
are pairwise non-proportional,
since $b$ does not divide
$p_n$ for these values of $n$.
Hence $\Card {\cal N''}\leq B(d,\eps )$.

To deal with $n\in{\cal N'''}$, we observe that by combining (7.5)
again with $z_1b^{t_n}+z_2b^{r_n}+z_3p_n=0$,
but now eliminating $p_n$, we obtain
$$
\biggl| {\xi z_3 +z_1 \over \xi z_3 - z_2} - {1\over b^{t_n}}
\biggr| < \biggl| {z_3 \over \xi z_3 - z_2} \biggr| \cdot
(b^{t_n})^{-1 -\eps}.
$$
In precisely the same manner as above, one obtains that
the pairs $(b^{t_n},1)$ lie in at most $B(d,\eps )$ one-dimensional subspaces.
Since these pairs are pairwise non-proportional, it follows that
$\Card {\cal N'''}\leq B(d,\eps )$.

By combining the above we obtain
$$
\Card{\cal N}=r\leq k_0+\Card{\cal N''}+\Card{\cal N'''}
\leq k_0+2B(d,\eps ).
$$
In view of (7.15), (7.19), this is smaller than
$(\eps^{-1})^{3+\eta}$ for $N$ sufficiently large.
This proves the claim, hence Theorem 2.1. \cqfd

\vskip 7mm

\centerline{\bf 8. Proof of Theorem 3.1}

\vskip 7mm

We closely follow Section 4 of \cite{BuPad}.
Assume without loss of generality that
$$
{b-1 \over b} < \xi < 1.
$$
Define the increasing sequence of positive integers $(n_j)_{j \ge 1}$
by $a_1 = \ldots = a_{n_1}$, $a_{n_1} \not = a_{n_1 + 1}$ and
$a_{n_j+1} = \ldots = a_{n_{j+1}}$, $a_{n_{j+1}} \not= a_{n_{j+1} + 1}$
for $j \ge 1$. Observe that
$$
\nbdc(n, \xi, b) = \max \{j : n_j \le n \}
$$
for $n \ge n_1$, and that $n_j \ge j$ for $j \ge 1$. Define
$$
\xi_j := \sum_{k = 1}^{n_j} \, {a_k \over b^k} +
\sum_{k = n_j + 1}^{+ \infty} \, {a_{n_j + 1} \over b^k}
= \sum_{k = 1}^{n_j} \, {a_k \over b^k} +
{a_{n_j + 1} \over b^{n_j} (b-1)} \cdot
$$
Then,
$$
\xi_j = {P_j(b) \over b^{n_j} (b-1)},
$$
where $P_j(X)$ is an integer polynomial of degree at most $n_j$
whose constant coefficient $a_{n_j+1} - a_{n_j}$ is not
divisible by $b$.
That is, $b$ does not divide $P_j(b)$.
We have
$$
|\xi - \xi_j| < {1 \over b^{n_{j + 1}}},
$$
and this can be rewritten as
$$
\biggl|(b-1) \xi - {P_j(b) \over b^{n_j} } \biggr|
< {b-1 \over b^{n_{j+1}}}.  \eqno (8.1)
$$
By Liouville's inequality,
$$
\biggl|(b-1) \xi - {P_j(b) \over b^{n_j} } \biggr|
\ge \bigl( 2 H\bigl((b-1) \xi \bigr) b^{n_j} \bigr)^{-d},
$$
so, if
$$
n_j \ge U := 1 + 3 H\bigl( (b-1) \xi \bigr), \eqno (8.2)
$$
then
$$
n_{j+1} \le 2 d n_j. \eqno (8.3)
$$

In what follows, constants implied by the Vinogradov symbols $\ll$, $\gg$
are absolute.
We need the following lemma.

\proclaim Lemma 8.1.
Let $0 < \eps \le 1$ and let $j_1$ denote the smallest $j$
such that $n_j \ge \max\{ U, 5/\eps \}$. Then
$$
\Card \{j : j\ge j_1, n_{j+1} / n_j \ge 1 + 2 \eps \}
\ll \log (6d)\eps^{-3}\log \big(\eps^{-1}\log (6d)\big).
$$

\proof For the integers $j$ into consideration, we have
$$
b^{n_j} > \max\bigl\{ 2H\bigl( (b-1) \xi \bigr), 2^{4/\eps} \bigr\}.
$$
Further, by (8.1), $n_j \ge U$, we get
$$
\biggl|(b-1) \xi - {P_j(b) \over b^{n_j} } \biggr|
< {b-1 \over (b^{n_j})^{1 + 2 \eps}}
\le {1 \over (b^{n_j})^{1 + \eps}}. \eqno (8.4)
$$
Moreover, for every prime $\ell$ dividing $b$,
$$
|b^{n_j}|_{\ell} \le \bigl(b^{n_j}\bigr)^{\log |b|_{\ell} / \log b}.
\eqno (8.5)
$$
Since
$$
1 + \eps + \sum_{\ell \mid b} \, {- \log |b|_{\ell} \over \log b}
= 2 + \eps,
$$
Corollary 5.2 applied to (8.4), (8.5) yields that for the integers $j$
into consideration the pairs $(P_j(b), b^{n_j})$
lie in
$$
\ll \log (6d)\eps^{-3}\log\big(\eps^{-1}\log (6d)\big)
$$
one-dimensional linear subspaces of $\Q^2$. But these
pairs are non-proportional since $b$ does not divide $P_j(b)$.
The lemma follows. \cqfd

Let $j_0$ be the smallest $j$ such that $n_j \ge U$.
Let $J$ be an integer with
$$
J > \max\bigl\{n_{j_0}^3, (4d)^6 \bigr\}. \eqno (8.6)
$$
Let $j_2$ be the largest integer with
$$
n_{j_2} \le 6 d J^{1/3}. \eqno (8.7)
$$
Then since $n_{j_2} \ge n_{j_0} \ge U$, we have
$$
n_{j_2} \ge {n_{j_2 + 1} \over 2 d} \ge 3 J^{1/3}. \eqno (8.8)
$$
Now choose
$$
\eps_1 := \Big( {\log (6d)\log J\over J}\Big)^{1/3} \eqno (8.9)
$$
and let $k$ be any positive integer and $\eps_2, \ldots , \eps_{k-1}$
any reals such that
$$
\eps_1<\eps_2 < \ldots < \eps_{k-1} < \eps_k :=1.
  \eqno (8.10)
$$
We infer from (8.8) that
$$
n_{j_2} \ge \max\{U, 5/\eps_h\}, \quad
\hbox{for $h=1, \ldots , k$}. \eqno (8.11)
$$
Let $\cS_0 = \{j_2, j_2+1, \ldots, J\}$ and,
for $h = 1, \ldots, k$, let $\cS_h$ denote the set of positive
integers $j$ such that $j_2 \le j \le J$
and $n_{j+1} \ge (1 + 2 \eps_h) n_j$.
Further, let $T_h$ be the cardinality of $\cS_h$ for $h=1, \ldots , k$.
Obviously, $\cS_0 \supset \cS_1 \supset \ldots \supset \cS_k$ and
$$
\cS_0 = (\cS_0 \setminus \cS_1) \cup (\cS_1 \setminus \cS_2)  \cup \ldots
\cup (\cS_{k-1} \setminus \cS_k) \cup \cS_k.
$$
Now,
$$
\eqalign{
{n_J \over n_{j_2}} & = {n_J \over n_{J-1}} \times
{n_{J-1} \over n_{J-2}}  \times \ldots \times
{n_{j_2+1} \over n_{j_2}}  \cr
& = \prod_{h=0}^{k-1} \biggl( \prod_{j \in \cS_h \setminus \cS_{h+1}}
\, {n_{j+1} \over n_j} \, \biggr) \,
\biggl( \prod_{j \in \cS_h}
\, {n_{j+1} \over n_j} \, \biggr) \cr
& \le (1 + 2 \eps_1)^J \, \prod_{h=1}^{k-1} \,
(1 + 2 \eps_{h+1})^{T_h - T_{h+1}}
\, (2d)^{T_k}, \cr}
$$
where in the last estimate we have used (8.11) and (8.3).
Taking logarithms, we get
$$
\eqalign{
\log (n_J/n_{j_2}) & \le 2 \eps_1 J + 2 \sum_{h=1}^{k-1} \eps_{h+1}
(T_h - T_{h+1}) + T_k \log (2d) \cr
& \le 2 \eps_1 J + 2 \eps_2 T_1 + 2 \sum_{h=2}^{k-1} (\eps_{h+1} - \eps_h)
T_h  - 2 T_k + T_k \log (2d). \cr}
$$
In view of (8.11), we can apply Lemma 8.1, and obtain that
$$
T_h \ll \log (6d)\eps_h^{-3} \, \log \big(\eps_h^{-1}\log (6d)\big)
$$
for $h = 1, \ldots , k$. This gives
$$
\eqalign{
\log (n_J/n_{j_2}) \ll \eps_1 J
& + \log (6d)\eps_2 \eps_1^{-3}\log\big(\eps_1^{-1}\log (6d)\big)
\cr
& + \log (6d)\, \sum_{h = 2}^{k-1} \, \eps_h^{-3}\log \big(\eps_h^{-1}\log (6d)\big)
\cdot (\eps_{h+1} - \eps_h) \cr
& +\big(\log (6d)\big)^2\log\log (6d).\cr}
$$
Now, let $k$ tend to infinity and $\max_{1 \le h \le k-1} \,
(\eps_{h+1} - \eps_h)$ tend to zero. Then the sum converges
to a Riemann integral,
and, after a short computation,
using that in view of (8.6), (8.9) we have $\eps_1^{-1}\gg d$,
we get
$$
\log (n_J/n_{j_2}) \ll \eps_1 J + \log (6d)\eps_1^{-2}\log (\eps_1^{-1}).
$$
By (8.6) and (8.7), we have $n_{j_2} \le J^{1/2} \le n_J^{1/2}$,
so $n_J/ n_{j_2} \ge n_J^{1/2}$. Inserting our choice (8.9)
for $\eps_1$ and using (8.6), we get
$$
\log n_J \ll J^{2/3}(\log J)^{1/3} \big(\log (6d)\big)^{1/3},
$$
i.e.,
$$
J\gg (\log n_J)^{3/2} (\log\log n_J)^{-1/2} \big(\log (6d)\big)^{-1/2}.
$$
This proves Theorem 3.1.
\qed

\vskip 5mm

\centerline{\bf 9. Final remarks}

\vskip 7mm

We deduce from Corollary 5.2 an improvement of an
extension due to Mahler \cite{Mah61}
of a theorem of Cugiani \cite{Cug58}, see \cite{BuCug}
for further references on the Cugiani--Mahler Theorem.

Let $S_1$, $S_2$ be finite, possibly empty sets of prime numbers,
put $S:=\{\infty\}\cup S_1\cup S_2$,
let $\xi\in\ovQ$ be an algebraic number,
let $\varepsilon >0$,
and let $f_p$ $(p\in S)$ be reals such that
$$
f_p\geq 0\ \hbox{for $p\in S$,}\ \ \sum_{p\in S} f_p=2.
$$

Let $\eps : \Z_{\ge 1} \to \R_{> 0}$ be a non-increasing function.
We consider the system of inequalities
$$
\left\{
\matrix{
|\xi -{x\over y}|\leq y^{-f_{\infty}-\eps(y)},\hfill\cr
|x|_p\leq y^{-f_p}\ (p\in S_1)\hfill\cr
|y|_p\leq y^{-f_p}\ (p\in S_2)\hfill
}
\right\}
\hbox{in $(x,y)\in\Z^2$ with $y>0$ and ${\rm gcd} (x, y) = 1$.}
\eqno (9.1)
$$

Arguing as in \cite{BuCug}, we get the following improvement 
of Theorem 1 on page 169 of \cite{Mah61}, that we state 
without proof.
For a positive integer $m$, we denote by $\exp_m$ the $m$th
iterate of the exponential function and by $\log_m$ the function
that coincides with the $m$th iterate of the logarithm function
on $[\exp_m 1, + \infty)$ and that takes the value $1$ on
$(-\infty, \exp_m 1]$.

\proclaim Theorem 9.1.
Keep the above notation. Let $m$ be a positive integer, and 
$c$ be a positive real number. Set
$$
\eps(y) = c \, (\log_{m+1} y)^{-1/3} \, (\log_{m+2} y),
\quad \hbox{for $y \ge 1$}.
$$
Let $(x_j/y_j)_{j \ge 1}$ be the sequence 
of reduced rational solutions of (9.1)
ordered such that $1 \le y_1 < y_2 < \ldots$ Then either the sequence
$(x_j / y_j)_{j \ge 1}$ is finite or
$$
\limsup_{j \to + \infty} \, {\log_m y_{j+1} \over \log_m y_j}
= + \infty.
$$

Theorem 9.1 improves upon Mahler's result, which deals only
with the case $m=1$ and involves the very slowly decreasing function
$y \mapsto (\log_3 y)^{-1/2}$.

Theorem 9.1 can be compared with Theorem 2 from \cite{BuCug} that
deals with products of linear forms and involves a function $\eps$
that depends on the cardinality of $S_1 \cup S_2$.
Note that Theorem 6.5.10 from 
Chapter 6 of the monograph of Bombieri and Gubler \cite{BoGu},
given without proof, 
deals also with products of linear forms, but the
function $\eps$ occurring there does not involve
the cardinality of $S_1 \cup S_2$.

We can then proceed exactly as Mahler did (\cite{Mah61}, Theorem 3, page 178)
to construct new explicit examples
of transcendental numbers. 

\proclaim Theorem 9.2.
Let $b \ge 2$ be an integer.
Let $\theta$ be a real number with $0 < \theta < 1$.
Let ${\bf n} = (n_j)_{j \ge 1}$ be an
increasing sequence of positive integers satisfying $n_1 \ge 3$
and
$$
n_{j+1} \ge \biggl( 1 +   {\log \log n_j
\over (\log n_j)^{1/3} }  \, \biggr) \, n_j, \quad (j \ge 1).
$$
Let $(a_j)_{j \ge 1}$ be a sequence of positive
integers prime to $b$ such that
$$
a_{j+1} \le b^{\theta (n_{j+1} - n_j)}, \quad j \ge 1.
$$
Then the real number
$$
\xi  = \sum_{j \ge 1} \, a_j \, b^{-n_j}
$$
is transcendental.

We omit the proof of Theorem 9.2, which follows from 
Theorem 9.1 with $m=1$.

It is of interest to note that Theorem 9.2 yields 
Corollary 3.2 only for $\eta > 3/4$. We would have 
obtained the same result by taking $k=1$
in (8.10). It 
is precisely the introduction of the
parameter $k$ there that allows us to get in Theorem 3.1
the exponent of $(\log n)$ equal to $3/2$ and not to $4/3$.

\vskip 5mm

\centerline{\bf APPENDIX}
\vskip 2mm
\centerline{\bf A quantitative two-dimensional Parametric Subspace
Theorem}

\vskip 7mm

We give a proof of the two-dimensional case of Proposition 4.1.
We keep the notation and assumptions from Section 4, except that we
assume $n=2$.
As before, $\K$ is an algebraic number field.
We recall the notation from Section 4, but now specialized to $n=2$.
Thus, ${\cal L}=( L_{iv}:\, v\in M_{\K},\, i=1,2)$ is a tuple of linear
forms satisfying
$$
\eqalignno{
&L_{iv}\in \K [X_1,X_2]\ \hbox{for $v\in\MK$, $i=1,2$,}&(A.1)
\cr
&L_{1v}=X_1,\, L_{2v}=X_2\ \hbox{for all but finitely many $v\in M_{\K}$,}
&(A.2)
\cr
&\det (L_{1v},L_{2v})=1\ \hbox{for } v\in M_{\K},&(A.3)
\cr
&\Card\Big(\bigcup_{v\in M_{\K}}\{ L_{1v},L_{2v}\}\Big)\leq r&(A.4)
}
$$
and
${\bf c}=(c_{iv}:\, v\in M_{\K},\, i=1,2)$ is a tuple of reals
satisfying
$$
\eqalignno{
&c_{1v}=c_{2v}=0\ \hbox{for all but finitely many $v\in M_{\K}$,}&(A.5)
\cr
&\sum_{v\in\MK}\sum_{i=1}^2 c_{iv}=0,&(A.6)
\cr
&\sum_{v\in\MK}\max (c_{1v},c_{2v})\leq 1.&(A.7)
}
$$
We define
$$
\HH =\HH (\LL ) :=
\prod_{v\in M_{\K}}\max_{1\leq i_1<i_2\leq s}\|\det (L_{i_1},L_{i_2})\|_v
\eqno (A.8)
$$
where we have written $\{ L_1\kdots L_s\}$ for
$\bigcup_{v\in M_{\K}}\{ L_{1v}, L_{2v}\}$.
Finally, for any finite extension $\E$ of $\K$
and any place $w\in\ME$ we define
$$
L_{iw}=L_{iv},\ c_{iw}=d(w|v)c_{iv}\quad\hbox{for } i=1,2,
\eqno (A.9)
$$
where $v$ is the place of $\MK$ lying below $w$ and $d(w|v)$ is given by (4.2).

The twisted height $H_{Q,\LL ,\cc}(\xx )$ of $\xx\in\ovQ^2$
is defined by taking any finite extension $\E$ of $\K$ such that $\xx\in\E^2$
and putting
$$
\HQLc (\xx ):=\prod_{w\in\ME}\max_{i=1,2}\| L_{iw}(\xx )\|_wQ^{-c_{iw}};
\eqno (A.10)
$$
this does not depend on the choice of $\E$.

\proclaim Proposition A.1.
Let $\LL =(L_{iv}:\, v\in\MK ,\, i=1,2)$
be a tuple of linear forms and
$\cc =(c_{iv}:\, v\in\MK ,\, i=1,2)$
a tuple of reals satisfying (A.1)--(A.7). Further, let $0<\delta \leq 1$.
\hfb
Then there are one-dimensional linear subspaces $T_1\kdots T_{t_2}$ of $\ovQ^2$,
all defined over $\K$, with
$$
t_2=t_2(r,\delta )=
2^{25}\delta^{-3}\log (2r)\log\big(\delta^{-1}\log (2r)\big)
\eqno (A.11)
$$
such that the following holds: for every real $Q$ with
$$
Q> \max \Big(\HH^{{2\over r(r-1)}},4^{1/\delta}\Big)\eqno (A.12)
$$
there is a subspace $T_i\in\{ T_1\kdots T_{t_2}\}$ such that
   
$$
\{ \xx\in\ovQ^2 :\, \HQLc (\xx )\leq Q^{-\delta}\}\subset T_i\, .
\eqno (A.13)
$$

The proof of Proposition A.1 is by combining some lemmata from \cite{EvSc02},
specialized to $n=2$. We keep the notation and assumptions from above.
By condition (A.4), there exists a `family' (unordered sequence with
possibly repetitions) of linear forms $\{ L_1\kdots L_r\}$
such that $L_{1v},L_{2v}$ belong to this family for every $v\in\MK$
and such that $L_1=X_1$, $L_2=X_2$. Now conditions (A.1)--(A.7) imply
the conditions (5.12)--(5.17) on p. 36 of \cite{EvSc02} with $n=2$.
These conditions are kept throughout \cite{EvSc02}
and so all arguments of \cite{EvSc02} from p. 36 onwards are applicable
in our situation.
Since in what follows the tuples $\LL$ and $\cc$ will be fixed and only
$Q$ will vary, we will write $H_Q$ for the twisted height $\HQLc$.

Let $Q$ be a real with $Q\geq 1$. We define the ``successive infima"
$\lone$, $\ltwo$ of $H_Q$ as follows: for $i=1,2$,
$\li$ is the infimum of all reals $\lambda >0$ such that
$\{ \xx\in\ovQ^2:\, H_Q(\xx )\leq\lambda\}$ contains at least $i$ linearly
independent points. Since we are working on 
the algebraic closure of $\Q$ and not on
a given number field, these infima need not be assumed by $H_Q$.

In \cite{EvSc02} (specialized to $n=2$), $\lone$, $\ltwo$ were defined to
be the successive infima of some sort of parallelepiped $\Pi (Q,\cc )$
defined over $\ovQ$,
and the lemmata in that paper were all formulated in terms
of these infima. 
However, according to \cite{EvSc02}, Corollary 7.4, p. 53, applied with $n=2$
and ${\bf A}=(Q^{c_{iv}},\, v\in\MK,\, i=1,2)$,
the successive infima of $\Pi (Q,\cc )$
are equal to the successive infima of $H_Q$ as defined above.

\proclaim Lemma A.2.
Let $\delta >0$, $Q\geq 1$.
\vskip 2mm\noindent
(i) ${1\over 2}\leq\lone\ltwo\leq 2$.
\vskip2mm\noindent
(ii) If there exists a non-zero $\xx\in\ovQ^2$ with $H_Q(\xx )\leq Q^{-\delta}$ then
$\lone\leq Q^{-\delta}$ and $\ltwo\geq {1\over 2}Q^{\delta}$.

\proof
Assertion (i) follows from \cite{EvSc02}, Corollary 7.6, p. 54.
Assertion (ii) is then obvious.
\qed

\proclaim Lemma A.3.
(Gap Principle). Let $\delta >0$, and let $Q_0$ be a real with
$Q_0> 4^{1/\delta}$. Then there is a unique, one-dimensional linear subspace $T$ of
$\ovQ^2$ with the following property: for every $Q$ with
$$
Q_0\leq Q< Q_0^{1+\delta /2}
$$
we have $\{ \xx\in\ovQ^2:\, H_Q(\xx )\leq Q^{-\delta}\}\subset T$.

\proof
Let $T$ be the linear subspace of $\ovQ^2$ spanned by all $\xx$ such that
$H_{Q_0}(\xx )\leq Q_0^{-\delta /2}$. If $T\not=({\bf 0})$ then by Lemma A.2
we have $\lambda_1(Q_0)\leq Q_0^{-\delta /2}$ and
$\lambda_2(Q_0)\geq {1\over 2}Q_0^{\delta /2}$,
which by our assumption on $Q_0$ is strictly larger than $\lambda_1(Q_0)$.
Hence $T$ has dimension at most $1$. So it suffices to prove that
if $\xx\in\ovQ^2$ and $Q$ are such that $Q_0\leq Q< Q_0^{1+\delta /2}$
and $H_Q(\xx )\leq Q^{-\delta}$, then $H_{Q_0}(\xx )\leq Q_0^{-\delta /2}$.

To prove this, choose a finite extension $\E$ of $\K$ such that $\xx\in\E^2$.
Notice that by (A.7), (A.9), (4.3) we have
$u:=\sum_{w\in\ME}\max (c_{1w},c_{2w})\leq 1$. For $w\in\ME$ we have
$$
\eqalign{
&\max\Big( {\| L_{1w}(\xx )\|_w\over Q_0^{c_{1w}}},
              {\| L_{2w}(\xx )\|_w\over Q_0^{c_{2w}}}\Big)\cr
&\quad
\leq \max\Big( {\| L_{1w}(\xx )\|_w\over Q^{c_{1w}}},
              {\| L_{2w}(\xx )\|_w\over Q^{c_{2w}}}\Big)\cdot
\Big({Q\over Q_0}\Big)^{\max (c_{1w},c_{2w})}\, .\cr
}
$$
So
$$
\eqalign{
H_{Q_0}(\xx ) &\leq H_Q(\xx )\Big({Q\over Q_0}\Big)^u
\cr
&\leq Q^{-\delta}\cdot {Q\over Q_0}\leq
Q_0^{-\delta}Q_0^{\delta /2}=Q_0^{-\delta/2}.
\cr
}
$$
$\quad$\qed

\proclaim Lemma A.4.
Let $\delta >0$ and let $A,B$ be reals with $4^{1/\delta}< A<B$.
\hfb
Then there are one-dimensional linear subspaces $T_1\kdots T_{t_3}$ of $\ovQ^2$,
with
$$
t_3\leq 1+\,{\log (\log B/\log A)\over \log (1+\delta /2)}
$$
such that for every $Q$ with $A\leq Q<B$ there is $T_i\in\{ T_1\kdots T_{t_3}\}$
with
$$
\{ \xx\in\ovQ^2:\, H_Q(\xx )\leq Q^{-\delta}\}\subset T_i.
$$

\proof
Let $k$ be the smallest integer with $A^{(1+\delta /2)^k}\geq B$.
Apply Lemma A.3 with $Q_0=A^{(1+\delta /2)^i}$ for $i=0\kdots k-1$.
\qed

We define the Euclidean height $H_2(\xx )$ for $\xx =(x_1\kdots x_m)\in\ovQ^m$
as follows. Choose any number field $\E$ such that $\xx\in\E^m$,
define
$$
\eqalign{
\|\xx\|_{w,2}:= 
\Big\{ \Big(\sum_{i=1}^m |x_i|_w^2\Big)^{1/2}\Big\}^{{[\E_w:\R ]\over [\E :\Q ]}} 
&\ \hbox{if $w$ is Archimedean,}\cr
\|\xx\|_{w,2}:=\max (\| x_1\|_w\kdots \| x_m\|_w) 
&\ \hbox{if $w$ is non-Archimedean,}\cr
}
$$
and put
$$
H_2(\xx ):=\prod_{w\in\ME}\|\xx \|_{w,2}.  
$$
This is independent of the choice of $\E$. For a polynomial $P$ with
coefficients in $\ovQ$, define $H_2(P):=H_2({\bf p})$,
where ${\bf p}$ is a vector consisting of the coefficients of $P$.

\proclaim Lemma A.5.
Let $\delta >0$. Consider the set of reals $Q$ such that
$$
\eqalignno{
&\hbox{there is $\xx_Q\in\ovQ^2\setminus\{ {\bf 0}\}$ 
with $H_Q(\xx_Q)\leq Q^{-\delta}$,}&(A.14)
\cr
&Q^{\delta}>(2\HH )^{6{r\choose 2}}.&(A.15)
}
$$
Then one of the following two alternatives is true:
\vskip2mm\noindent
(i) For all $Q$ under consideration we have
$H_2(\xx_Q)> Q^{\delta /3{r\choose 2}}$;
\vskip2mm\noindent
(ii) There is a single one-dimensional linear subspace $T_0$ of $\ovQ^2$
such that for all $Q$ under consideration we have $\xx_Q\in T_0$.

\proof
This is \cite{EvSc02}, p.80, Lemma 12.4 with $n=2$.
Condition (A.14) and Lemma A.2 imply $\lone\leq Q^{-\delta}$ which is
condition (12.37) of Lemma 12.4 of \cite{EvSc02} with $n=2$.
Further, the quantity $R$ in that lemma is $\leq {r\choose 2}$
(see \cite{EvSc02}, p.75, Lemma 12.1).
\qed

Let $m$ be a positive integer and $\rr =(r_1\kdots r_m)$ a tuple of positive
integers. We say that a polynomial is multihomogeneous of degree $\rr$
in the blocks of variables $\X_1=(X_{11},X_{12})\kdots \X_m=(X_{m1},X_{m2})$
if it can be expressed as a linear combination of monomials
$$
\prod_{h=1}^m\prod_{k=1}^2 X_{hk}^{i_{hk}}\ \
\hbox{with $i_{h1}+i_{h2}=r_h$ for $h=1\kdots m$.}
$$
(Below, $h$ will always index the block).
Given points
$\xx_h=(x_{h1},x_{h2})$ $(h=1\kdots m)$ 
and a polynomial $P$ which is multihomogeneous in $\X_1\kdots \X_m$, 
we write $P(\xx_1\kdots\xx_m)$
for the value obtained by substituting $x_{hk}$
for $X_{hk}$ ($h=1\kdots m$, $k=1,2$).

The {\it index} of a polynomial $P$ multihomogeneous in $\X_1\kdots\X_m$
with respect to points $\xx_1\kdots\xx_m$
and to a tuple of positive integers $\rr =(r_1\kdots r_m)$,
denoted by
$$
\indP ,
$$
is defined to be the smallest
real $\sigma$ with the following property: there is a tuple of non-negative
integers ${\bf i}=(i_{hk}:\, h=1\kdots m,\, k=1,2)$ such that
$$
\eqalign{
&\Biggl(\prod_{h=1}^m\prod_{k=1}^2
\Big({\partial\over \partial X_{hk}}\Big)^{i_{hk}}P\Biggr) (\xx_1\kdots\xx_m)\not= 0;
\cr
&\sum_{h=1}^m {i_{h1}+i_{h2}\over r_h}=\sigma\, .
}
$$
For a field ${\bf F}$ and a tuple of positive integers $\rr =(r_1\kdots r_m)$,
We denote by ${\bf F} [\rr ]$ the set of polynomials
with coefficients in ${\bf F}$
which are multihomogeneous of degree $\rr$ in $\X_1\kdots\X_m$.

We define the constant $C(\K ):= |D_{\K}|^{1/[\K:\Q ]}$,
where $D_{\K}$ denotes the discriminant of $\K$.
In fact, the precise value of $C(\K )$ is not of importance. 

\proclaim Lemma A.6.
Suppose that $0<\delta\leq 1$, let $\theta$ be a real with
$$
0<\theta\leq {\delta\over 80},
\eqno (A.16)
$$
$m$ an integer with
$$
m> 4\theta^{-2}\log (2r)
\eqno (A.17)
$$
and $\rr =(r_1\kdots r_m)$ a tuple of positive integers, and
put $q:= r_1+\cdots +r_m$.
\hfill\break  
Suppose that there exist
positive reals $Q_1\kdots Q_m$ and non-zero points $\xx_1\kdots\xx_m$
in $\ovQ^2$ such that
$$
\eqalignno{
&r_1\log Q_1\leq r_h\log Q_h\leq (1+\theta )r_1\log Q_1\ \ (h=1\kdots m),&(A.18)
\cr
&H_{Q_h}(\xx_h)\leq Q_h^{-\delta}\ \ (h=1\kdots m),&(A.19)
\cr
&Q_h^{\delta}> C(\K )^{5/4q}\cdot 2^{50}\HH^5 \theta^{-5/2}.&(A.20)
}
$$
Then there is a non-zero polynomial $P\in K[\rr ]$ such that
$$
\eqalignno{
&\indP \geq m\theta ,&(A.21)
\cr
&H_2(P)\leq C(\K )^{1/2}\cdot 2^{3m} (12\HH )^q.&(A.22)
}
$$

\proof
This is \cite{EvSc02}, Lemma 15.1, p. 89, with $n=2$. The space $V_{[h]}(Q_h)$
in that lemma is in our situation precisely the space spanned by $\xx_h$
for $h=1\kdots m$. Inequality (A.17) comes from (14.7) on \cite{EvSc02}, p. 83;
later it is assumed that $s=r$ (see (14.10) on \cite{EvSc02}, p.85).
Inequality (A.22) comes from the inequality at the bottom of p. 87 of
\cite{EvSc02}. The construction of the polynomial $P$ is by means of a
now standard argument, based on the Bombieri--Vaaler Siegel's Lemma.
\qed

\proclaim Lemma A.7.
(Roth's Lemma) Let $0<\theta \leq 1$. Let $m$ be an integer with $m\geq 2$
and $\rr =(r_1\kdots r_m)$ a tuple of positive integers such that
$$
{r_h\over r_{h+1}}\geq {2m^2\over\theta}\ \ (h=1\kdots m-1).
\eqno (A.23)
$$
Further, 
let $P$ be a non-zero polynomial in $\ovQ [\rr ]$ 
and $\xx_1\kdots\xx_m$
non-zero points in $\ovQ^2$ such that
$$
H_2(\xx_h)^{r_h}\geq \big( e^qH_2(P)\big)^{(3m^2/\theta)^m}\ \ (h=1\kdots m)
\eqno (A.24)
$$
where $e=2.7182\ldots$, $q=r_1+\cdots +r_m$. Then
$$
\indP <m\theta\, .\eqno (A.25)
$$

\proof
This is the case $n=2$ of \cite{Ev96}, Lemma 24. It is an immediate
consequence of \cite{Ev95}, Theorem 3.
\qed

We keep our assumption $0<\delta\leq 1$ and define the integer
$$
m:= 1+[ 25600\cdot\delta^{-2}\log (2r)].
\eqno (A.26)
$$
Put
$$
C:= (36\HH )^{m\cdot (240m^2/\delta )^m\cdot 3{r\choose 2}/\delta}. 
\eqno (A.27)
$$
Denote by $\SS$ the set of reals $Q$ such that
$$
Q\geq C;\ \ \hbox{there is $\xx\in\ovQ^2\setminus\{ {\bf 0}\}$
with $H_Q(\xx )\leq Q^{-\delta}$.}
\eqno (A.28)
$$

\proclaim Lemma A.8.
One of the following two alternatives is true:
\vskip2mm\noindent
(i) There is a single, one-dimensional linear subspace $T_0$ of $\ovQ^2$
such that for every $Q\in\SS$ we have 
$\{ \xx\in\ovQ^2:\, H_Q(\xx)\leq Q^{-\delta}\}\subset T_0$;
\vskip2mm\noindent
(ii) There are reals $Q_1\kdots Q_{m-1}$ with $C\leq Q_1<\cdots <Q_{m-1}$
such that
$$
\SS \subset \bigcup_{h=1}^{m-1} [Q_h,Q_h^{162m^2/\delta}].
\eqno (A.29)
$$

\proof
We assume that neither of the alternatives (i) or (ii) is true.
From this assumption, we will deduce that there are a tuple of positive
integers $\rr =(r_1\kdots r_m)$, a non-zero polynomial $P\in \K [\rr ]$,
and non-zero points $\xx_1\kdots\xx_m\in\ovQ^2$, satisfying both
(A.21) and (A.25). This is obviously impossible.

By our assumption, there are reals $Q_1\kdots Q_m\in\SS$ with 
$$
{\log Q_{h+1}\over \log Q_h}\geq {162m^2\over\delta}\ \
(h=1\kdots m-1),
\eqno (A.30)
$$
and non-zero points $\xx_1\kdots\xx_m\in\ovQ^2$ with
$$
H_{Q_h}(\xx_h)\leq Q_h^{-\delta}\ \ (h=1\kdots m).\eqno (A.31)
$$
Put $\theta :=\delta /80$. First choose a positive integer $s_1$ such that
$\theta s_1\log Q_1 >\log Q_h$ for $h=2\kdots m$. Then there are integers
$s_2\kdots s_m$ such that
$$
s_1\log Q_1\leq s_h\log Q_h\leq (1+\theta )s_1\log Q_1\ \ (h=1\kdots m).
$$
Now take $r_h:= ts_h$ $(h=1\kdots m)$, $\rr =(r_1\kdots r_m)$, where $t$
is a positive integer, chosen large enough such that the right-hand side
of (A.20) is smaller than $C^{\delta}$ and the right-hand side of (A.22)
is smaller than $(13\HH )^q$, where $q=r_1+\cdots +r_m$. Then conditions
(A.18)--(A.20) of Lemma A.6 are satisfied, hence there exists a non-zero
polynomial $P\in \K [\rr ]$ such that (A.21), (A.22) are satisfied.
So we have in fact
$$
H_2(P)\leq (13\HH )^q .\eqno (A.32)
$$

We now show that $P$, $\rr$, $\xx_1\kdots\xx_m$ satisfy conditions (A.23),
(A.24) of Lemma A.7. Then it follows that (A.25) holds, and we arrive at the
contradiction we wanted.

In view of (A.30), (A.18) and $\theta ={\delta \over 80}\leq {1\over 80}$
we have
$$
\eqalign{
{r_h\over r_{h+1}} &={r_h\log Q_h\over r_{h+1}\log Q_{h+1}}
                             \cdot {\log Q_{h+1}\over\log Q_h}
\cr
&\ge {1\over 1+\theta}\cdot {162m^2\over \delta} 
                        \geq {160m^2\over \delta}={2m^2\over\theta}
\cr
}
$$
for $h=1\kdots m-1$, which is (A.23).

Our reals $Q\in\SS$ satisfy conditions (A.14), (A.15) of Lemma A.5.
Since we assumed that alternative (i) of Lemma A.8 is false,
alternative (ii) of Lemma A.5 must be false. So alternative (i) of that
lemma must be true. This implies in particular, that
$$
H_2(\xx_h) > Q_h^{\delta /2{r\choose 2}}\ \ (h=1\kdots m).
$$
By combining this with (A.18), (A.32), this implies
$$
\eqalign{
H_2(\xx_h)^{r_h} &\geq (Q_h^{r_h})^{\delta/3{r\choose 2}}\geq
(Q_1^{r_1})^{\delta/3{r\choose 2}}\geq C^{r_1\delta/3{r\choose 2}}
\cr
&\geq (36\HH )^{mr_1(3m^2/\theta )^m}\geq (e^qH_2(P))^{(3m^2/\theta )^m}
}
$$
for $h=1\kdots m$, which is (A.25). This completes our proof.
\qed

\noindent
{\it Proof of Proposition A.1.}
First suppose that alternative (ii) of Lemma A.8 is true.
By applying Lemma A.4 with $A=Q_h$, $B=Q_h^{162m^2/\delta}$ for $h=1\kdots m-1$
we conclude the following:
\vskip2mm
There are one-dimensional linear subspaces $T_1\kdots T_{t_4}$ of
$\ovQ^2$, with
$$
t_4\leq (m-1)\Big\{ 1+{\log (162m^2/\delta )\over\log (1+\delta /2)}\Big\}
\leq 5\delta^{-1}m\log (162m^2/\delta )
$$
such that for every $Q$ with
$$
Q\geq C :=(36\HH )^{(240m^2/\delta )^m\cdot 3m{r\choose 2}/\delta}
$$
there is $T_i\in\{ T_1\kdots T_{t_4}\}$ with
$\{ \xx\in\ovQ^2:\, H_Q(\xx )\leq Q^{-\delta}\}\subset T_i$.
This holds true trivially also if alternative (i) of Lemma A.8 is true;
so it holds true in all cases.

It remains to consider those values $Q$ with
$$
\max \big(\HH^{1/{r\choose 2}},4^{1/\delta}\big)=: C'< Q<C. 
\eqno (A.33)
$$
Notice that $C'\geq (36\HH )^{1/12{r\choose 2}}$. Hence by Lemma A.4,
there are one-dimensional linear subspaces
$T_1'\kdots T_{t_5}'$ of $\ovQ^2$, with
$$
\eqalign{
t_5 &\leq 1+{\log (\log C/\log C')\over \log (1+\delta /2)}\cr
&\leq 5\delta^{-1}\Big( m\log (240m^2/\delta \Big)+
\log \Big(3m{r\choose 2}\Big)+\log \Big(12{r\choose 2}\Big)\cr
&\leq 6\delta^{-1}m\log (240m^2/\delta )
}
$$
such that for every $Q$ with (A.33) there is $T_i'\in\{ T_1'\kdots T_{t_5}'\}$
with\hfb $\{\xx\in\ovQ^2:\, H_Q(\xx )\leq Q^{-\delta}\}\subset T_i'$.

Collecting the above, we get that there are one-dimensional linear subspaces
\hfill\break
$T_1\kdots T_{t_2}$ of $\ovQ^2$, with
$$
t_2 \leq t_4+t_5\leq 11\delta^{-1}m\log (240m^2/\delta )
\leq 33\delta^{-1}m\log m
$$
such that for every $Q> C'$ there is $T_i\in\{ T_1\kdots T_{t_2}\}$
with
\hfb
$\{ \xx\in\ovQ^2:\, H_Q(\xx )\leq Q^{-\delta}\}\subset T_i$.
Substituting (A.26) for $m$ we obtain
$$
\eqalign{
t_2 &\leq 33\delta^{-1}\cdot 25601\delta^{-2}\log (2r)
\log(25601\delta^{-2}\log (2r))\cr
&< 2^{25}\delta^{-3}\log (2r)\log \big(\delta^{-1}\log (2r)\big)
}
$$
which is the right-hand side of (A.11).

To finish the proof of Proposition A.1, it remains to show that the spaces
$T_1\kdots T_{t_2}$ are defined over $\K$.
Let $Q$ be any real $\geq 1$. Suppose that there are non-zero
vectors $\xx\in\ovQ^2$ with $H_Q(\xx )\leq Q^{-\delta}$,
and that these vectors span a one-dimensional linear subspace $T$ of $\ovQ^2$.
According to \cite{EvSc02}, Lemma 4.1, p.32, we have
for any $K$-automorphism $\sigma$ of $\ovQ$ that
$H_Q(\sigma (\xx ))=H_Q(\xx )$,
where $\sigma (\xx )$ is obtained by applying $\sigma$ to the coordinates
of $\xx$; hence $\sigma (\xx )\in T$.
This implies that $T$ is defined over $\K$.
\qed

\vskip 9mm


\centerline{\bf References}

\vskip 7mm

\beginthebibliography{999}

\bibitem{AdBu07}
B. Adamczewski and Y. Bugeaud,
{\it On the complexity of algebraic numbers I. Expansions in integer bases},
Ann. of Math. 165 (2007), 547--565.

\bibitem{AdBuLu}
B. Adamczewski, Y. Bugeaud et  F. Luca,
{\it Sur la complexit\'e des nombres alg\'ebriques},
C. R. Acad. Sci. Paris 339 (2004), 11--14.

\bibitem{AlbTh}
J. Albert,
Propri\'et\'es combinatoires et arithm\'etiques de
certaines suites automatiques et substitutives.
Th\`ese de doctorat, Universit\'e Paris XI, 2006.

\bibitem{All00}
J.-P. Allouche,
{\it Nouveaux r\'esultats de transcendance de r\'eels \`a
d\'eveloppements non al\'eatoire},
Gaz. Math. 84 (2000), 19--34.

\bibitem{BBCP}
D. H. Bailey, J. M. Borwein, R. E. Crandall and C. Pomerance, 
{\it On the binary expansions of algebraic numbers},
J. Th\'eor. Nombres Bordeaux 16 (2004), 487--518.

\bibitem{BoGu}
E. Bombieri and W. Gubler,
Heights in Diophantine Geometry.
New mathematical monographs 4, Cambridge University Press, 2006.


\bibitem{Bor50}
\'E. Borel,
{\it Sur les chiffres d\'ecimaux de $\sqrt{2}$ et divers
probl\`emes de probabilit\'es en cha\^\i ne},
C.~ R.~ Acad.~ Sci.~ Paris 230 (1950), 591--593.

\bibitem{BuPad}
Y. Bugeaud,
{\it On the $b$-ary expansion of an algebraic number},
Rend. Sem. Math. Univ. Padova. To appear.

\bibitem{BuCug}
Y. Bugeaud,
{\it Extensions of the Cugiani--Mahler Theorem}.
Preprint.

\bibitem{Cha33}
D. G. Champernowne,
{\it The construction of decimals normal in the scale of ten},
J. London Math. Soc. 8 (1933), 254--260.

\bibitem{Cug58}
M. Cugiani,
{\it Sull'approssimazione di numeri algebrici mediante razionali},
Collectanea Mathematica,
Pubblicazioni dell'Istituto di matematica dell'Universit\`a di Milano 169,
Ed. C. Tanburini, Milano,
pagg. 5 (1958).

\bibitem{Ev95}
J.-H. Evertse,
{\it An explicit version of Faltings' Product Theorem and an improvement
of Roth's Lemma}, Acta Arith. 73 (1995), 215--248.

\bibitem{Ev96}
J.-H. Evertse,
{\it An improvement of the quantitative Subspace Theorem},
Compos. Math. 101 (1996), 225--311.


\bibitem{EvSc02}
J.-H. Evertse and H.P. Schlickewei,
{\it A quantitative version of the Absolute Subspace Theorem},
J. reine angew. Math. 548 (2002), 21--127.

\bibitem{FeMa97}
S. Ferenczi and Ch. Mauduit,
{\it Transcendence of numbers with a low complexity expansion},
J. Number Theory 67 (1997), 146--161.

\bibitem{Lo99}
H. Locher,
{\it On the number of good approximations of algebraic
numbers by algebraic numbers of bounded degree},
Acta Arith. 89 (1999), 97--122.

\bibitem{Mah61}
K. Mahler,
Lectures on Diophantine approximation,
Part 1: $g$-adic numbers and Roth's theorem,
University of Notre Dame, Ann Arbor, 1961.

\bibitem{MoHe38}
M. Morse and G. A. Hedlund,
{\it Symbolic dynamics},
Amer. J. Math. 60 (1938), 815--866.

\bibitem{MoHe40}
M. Morse and G. A. Hedlund,
{\it Symbolic dynamics II},
Amer. J. Math. 62 (1940), 1--42.

\bibitem{Pan84}
J.-J. Pansiot,
{\it Bornes inf\'erieures sur la complexit\'e des
facteurs des mots infinis engendr\'es par morphismes it\'er\'es},
STACS 84 (Paris, 1984),  230--240,
Lecture Notes in Comput. Sci., 166, Springer, Berlin, 1984.

\bibitem{Rid57}
D. Ridout,
{\it Rational approximations to algebraic numbers},
Mathematika 4 (1957), 125--131.

\bibitem{Ri06}
T. Rivoal,
{\it On the bits counting function of real numbers},
J. Austral. Math. Soc. To appear.

\bibitem{Ro55}
K. F. Roth,
{\it Rational approximations to algebraic numbers},
Mathematika 2 (1955), 1--20; corrigendum, 168.

\bibitem{SchmLN}
W. M. Schmidt,
{Diophantine Approximation}.
Lecture Notes in Mathematics 785, Springer, 1980.

\bibitem{Schn}
{T. Schneider},
Einf\"uhrung in die transzendenten Zahlen.
Springer--Verlag, Berlin--G\"ottingen--Heidelberg, 1957.

\endthebibliography

\vskip 10mm

\noindent{Yann Bugeaud}

\noindent{Universit\'e Louis Pasteur}

\noindent{U. F. R. de math\'ematiques}

\noindent{7, rue Ren\'e Descartes}

\noindent{67084 STRASBOURG Cedex (FRANCE)}

\vskip2mm

\noindent{{\tt bugeaud@math.u-strasbg.fr}}

\vskip 5mm

\noindent{Jan-Hendrik Evertse}

\noindent{Universiteit Leiden}

\noindent{Mathematisch Instituut}

\noindent{Postbus 9512}

\noindent{2300 RA LEIDEN (THE NETHERLANDS)}

\vskip2mm

\noindent{{\tt evertse@math.leidenuniv.nl}}

\goodbreak

\bye